\documentclass[a4paper,11pt,twoside]{article}

\usepackage{float}
\usepackage[utf8]{inputenc}
\usepackage[T1]{fontenc}
\usepackage[english]{babel}

\usepackage{charter}

\usepackage[dvipsnames]{xcolor}
\usepackage{verbatim}
\usepackage{listings}
\usepackage{multirow}

\usepackage[top=2.5cm, bottom=4.0cm, left=1.8cm, right=1.8cm]{geometry}
\usepackage{amsmath}
\usepackage{amsthm}	
\usepackage{amsfonts}	
\usepackage{amssymb	
	,bbm
	,units 
}
\usepackage{enumerate}
\usepackage[nobysame
,alphabetic
,initials
]{amsrefs}

\usepackage{hyperref}

\usepackage{algorithm}
\usepackage[noend]{algpseudocode}

\makeatletter
\renewcommand{\ALG@beginalgorithmic}{\footnotesize}
\makeatother

\usepackage{wrapfig}
\usepackage{graphicx}

\numberwithin{equation}{section}
\numberwithin{figure}{section}
\numberwithin{algorithm}{section}
\numberwithin{table}{section}
 \usepackage[nodayofweek]{datetime}

\renewcommand*{\thefootnote}{\fnsymbol{footnote}}

\title{An unbiased It\^{o} type stochastic representation for transport PDEs: A Toy Example}

\author{Greig Smith}
\author{
 	\normalsize Gon\c calo dos Reis \footnote{G. dos Reis acknowledges support from the \emph{Funda{\c c}$\tilde{\text{a}}$o para a Ci$\hat{e}$ncia e a Tecnologia} (Portuguese Foundation for Science and Technology) through the project UID/MAT/00297/2019 (Centro de Matem\'atica e Aplica\c c$\tilde{\text{o}}$es CMA/FCT/UNL).} \\[8pt]
         \small  University of Edinburgh\\ 
         \small  and \\
	\small  Centro de Matem\'atica e Aplica\c c$\tilde{\text{o}}$es\\
        \small  G.dosReis@ed.ac.uk
 \and
          \normalsize Greig Smith \footnote{G. Smith was supported by The Maxwell Institute Graduate School in Analysis and its
          	Applications, a Centre for Doctoral Training funded by the UK Engineering and Physical
          	Sciences Research Council (grant EP/L016508/01), the Scottish Funding Council, the University of Edinburgh and Heriot-Watt
          	University.}  \\[8pt]
          \small  University of Edinburgh\\
          \small Maxwell Institute for Mathematical Sciences\\
          \small School of Mathematics\\
          \small  G.Smith-13@sms.ed.ac.uk
}

\renewcommand{\thefootnote}{\arabic{footnote}}

 \date{\currenttime, \ddmmyyyydate\today
      }

\theoremstyle{plain}
\newtheorem{theorem}{Theorem}[section]
\newtheorem{lemma}[theorem]{Lemma}
\newtheorem{proposition}[theorem]{Proposition}

\newtheorem{remark}[theorem]{Remark}

\newtheorem{assumption}[theorem]{Assumption}

\newcommand{\bE}{\mathbb{E}}

\newcommand{\bI}{\mathbb{I}}

\newcommand{\bN}{\mathbb{N}}
\newcommand{\bP}{\mathbb{P}}

\newcommand{\bR}{\mathbb{R}}
\newcommand{\bS}{\mathbb{S}}

\newcommand{\cF}{\mathcal{F}}

\newcommand{\cL}{\mathcal{L}}

\newcommand{\cW}{\mathcal{W}}

\newcommand{\dd}{\mathrm{d}}

\newcommand{\1}{\mathbbm{1}}

\begin{document}

\selectlanguage{english}

\maketitle

\renewcommand*{\thefootnote}{\arabic{footnote}}

\begin{abstract}
We propose a stochastic representation for a simple class of transport PDEs based on It\^o representations. We detail an algorithm using an estimator stemming for the representation that, unlike regularization by noise estimators, is unbiased. We rely on recent developments on branching diffusions, regime switching processes and their representations of PDEs. 

There is a loose relation between our technique and regularization by noise, but contrary to the latter, we add a perturbation and immediately its correction. The method is only possible through a judicious choice of the diffusion coefficient $\sigma$.  A key feature is that our approach does not rely on the smallness of $\sigma$, in fact, our $\sigma$ is strictly bounded from below which is in stark contrast with standard perturbation techniques. This is critical for extending this method to non-toy PDEs which have nonlinear terms in the first derivative where the usual perturbation technique breaks down. 
 
The examples presented show the algorithm outperforming alternative approaches. Moreover, the examples point toward a potential algorithm for the fully nonlinear case where the method of characteristics break down.
\end{abstract}



\footnotesize

\normalsize


\section{Introduction}

Stochastic techniques to solve PDEs have become increasing popular in recent times with advances in computing power and numerical techniques allowing for solutions of PDEs to be calculated to high precision. Advances in BSDEs (Backward Stochastic Differential Equations) and so-called branching diffusions also allow one to tackle nonlinear PDEs (see \cite{BernalEtAl2017} and references therein). Stochastic representations for PDEs are useful as they give access to probabilistic Monte Carlo methods, in turn yielding strong  numerical gains over deterministic based solvers, especially in high dimensional problems, see \cites{han2017overcoming,fahim2011probabilistic,BernalEtAl2017}. Unlike their deterministic counterparts, stochastic based PDE solvers are less prone to the curse of dimensionality. In \cite{BernalEtAl2017} the authors used hybrid Monte Carlo \& PDE solvers to split the domain of the non-linear PDEs into multiple (independent) parts which allowed one to achieve perfect parallelization drastically reducing the time taken to numerically approximate such equations; a general discussion on such techniques is given there. 

In this work we focus on transport PDEs. One of the main limitations when using It\^o based stochastic techniques to represent PDEs is the requirement that the PDE is of second order in space (i.e.~a ``Laplacian'' must be present). Thus PDEs with only one spatial and one time derivative (transport PDEs) have been, until now beyond the scope of stochastic techniques. An idea to navigate around this is to perturb the PDE by a ``small'' Laplacian, then one can use stochastic techniques on the perturbed PDE. Although this does provide a way to approximate the solution, it is very dependent on the perturbation being small enough so that the solution of the perturbed PDE is close to the first order PDE. Of course introducing a perturbation will lead to an error (bias) in the estimation, but more problematic is that the inverse of the perturbation coefficient will appear in the nonlinearities containing derivatives, thus the small perturbation makes the numerical scheme unstable. We discuss this point further in Section \ref{sec:examples}. Let us note that stochastic representations are only important for transport PDEs with nonlinearities in the derivative of the solution, see Remark \ref{Rem:Non stochastic representation}.

A string of related literature based on numerical approximations via branching processes has re-emerged due to to recent developments. We do not carry out a review of these developments here but refer to \cite{BernalEtAl2017} for a review on the state of the art. Branching algorithms offer a useful approach to solve non-linear PDEs and also for unbiased simulation of SDEs (see \cites{LabordereEtAl2016,DoumbiaEtAl2017}). However, in order to apply Monte Carlo methods one requires estimators to be square-integrable and of finite computational complexity. For square integrability several works have fine tuned previous results to allow for increasing general cases: \cite{LabordereEtAl2015} introduced a control variate on the final step, which allowed for an unbiased simulation of an SDE with constant diffusion; later, \cite{LabordereEtAl2016} changed the time stepping scheme from an Exponential to a Gamma random variable, this allowed for the simulation of semilinear PDEs; most recently, \cite{DoumbiaEtAl2017} used antithetic variables as well as control variates to obtain an unbiased algorithm for an SDE with non constant diffusion. 

The material we present requires all of the above mentioned improvements along with new ideas in order to ensure the estimator to be square-integrable. Taking the long view, we believe these techniques to be crucial in extending this type of stochastic representations to the fully non-linear case. The second order parabolic fully nonlinear case has been considered in \cite{Warin2017} and \cite{Warin2018}, but the theoretical basis for that case is to the best of our knowledge open. There are also several works looking at branching style algorithms but to tackle different types of PDEs, see \cite{CuchieroTeichmann2017}, \cite{AgarwalClaisse2017} and \cite{LabordereTouzi2018} for further results.
\medskip

The \emph{contributions} of this paper are two-fold. Firstly we show how one can take the ideas of branching diffusions and regime switching to construct an unbiased stochastic representation for transport PDE. To the best of our knowledge this is the first result of its kind. Secondly, we improve upon the techniques currently presented in the literature \cites{LabordereEtAl2016,DoumbiaEtAl2017,LabordereEtAl2015} in order to show our representation is square integrable and of finite computational complexity and thus can be used in Monte Carlo simulation. For better readability we also provide a heuristic description of our ideas.

From a \emph{methodological} point of view, the approach in this paper is related to the regime switching algorithms presented in \cite{DoumbiaEtAl2017} and \cite{LabordereEtAl2015}, where one adds and subtracts terms in the PDE to change the ``driving SDE'' defined by the Dynkin operator. Such algorithms were inspired by branching diffusion algorithms as developed in \cite{RasulovRaimovaMascagni2010} and \cite{Labordere2012}, although there is also a connection with parametrix approach (see \cite{AnderssonHiga2017}) where measure changes are used with corresponding weights to yield an unbiased representation. Here we add and subtract the second order derivative, which leaves us with a nonlinear PDE that can then be solved using regime switching (essentially we perturb the PDE then correct for the perturbation). Crucially this does not require $\sigma$ to be small. Although the transport PDE we consider is simple, one of the main challenges is to keep the representation square integrable, which comes from the added second order term. The general case will be addressed in future work, nonetheless we give numerical examples showing that the general case is within (numerical) reach.

\begin{remark}
	\label{Rem:Non stochastic representation}
	Basic first order PDEs can easily be made to have a stochastic like representation using branching type arguments, for example a PDE of the type,
	\begin{align*}
	\partial_{t} u(t,x) + b(t,x)\partial_{x}u(t,x) + u(t,x)^{2}=0, \quad u(T,x)=g(x) \, .
	\end{align*}
	It is possible to write the solution to this as,
	\begin{align*}
	u(t,x)=g(X_{T})+ \int_{t}^{T} u(s,X_{s})^{2} \dd s \, ,
	\end{align*}
	where $X$ is the deterministic process satisfying the ODE $\dd X_{s} = b(t,X_{s}) \dd s$, $X_{t}=x$. Introducing random times into the solution of $u$ as is done in standard branching we can obtain a solution to $u$ as the expected product of particles at time $T$. A similar argument can also be made for nonlinear ODEs. 
	
	What is crucial here though is that this argument only holds when we do not have nonlinearities in the first derivative of the process, since we require Malliavin integration by parts tricks to deal with those. This is also the case when we want to apply the unbiased trick to $b$.
\end{remark}

This work is organized as follows. In Section \ref{Sec:Representation} we present our notation, the problem and give a heuristic description of our ideas. In Section \ref{Sec:Main Results} we present and prove our main results. Following that in Section \ref{Sec:Toward the general case} we discuss the open problems left by this work. Finally Section \ref{sec:examples} illustrates numerically our findings to show our method is indeed unbiased. Moreover, we show the capability of our method to tackle problems in the nonlinear setting where the perturbation technique performs poorly.


\section{Regime Switching Diffusion Representation in general}
\label{Sec:Representation}
\subsection{Notation and recap of stochastic representations}
Following the standard notation in stochastic analysis let $C_{b}^{1,n}([0,T]\times \bR^{d},\bR)$ be the set of functions $v : [0,T]\times \bR^{d} \rightarrow \bR$ with one bounded time derivative and order $n$ bounded spatial derivatives. Further, let $d \ge 1$ and $W$ be a $d$-dimensional Brownian motion, defined on the probability space $(\Omega, \cF, \bP, (\cF_{t})_{t \ge 0})$, with $\cF_{t}$ the filtration of a multidimensional Brownian motion augmented with the null sets (satisfying the usual conditions). 

Consider a multidimensional stochastic differential equation (SDE) $X$ starting at time point $t$, $0 \le t \le T$ of the form,
\begin{align*}
\dd X_{s}= b(s,X_{s}) \dd s + \sigma(s, X_{s}) \dd W_{s} \, , \quad \text{for } s \in [t,T] ~ \text{and } X_{t}=x \, , 
\end{align*}
where the drift $b: [0,T]\times \bR^{d} \rightarrow \bR^{d}$ and diffusion $\sigma:[0,T]\times \bR^{d} \rightarrow \bS^{d}$ satisfy the usual Lipschitz conditions so that the above SDE has a unique strong solution, $\bS^{d}$ denotes the set of $d$-by-$d$ dimensional real valued matrices. 

We associate with the SDE the infinitesimal generator $\cL$, which when applied to any function $\phi \in C_{b}^{1,2}([0,T]\times \bR^{d},\bR)$ in the domain of $\cL$ is,
\begin{align*}
(\cL \phi)(t,x)=b(t,x) \cdot D \phi(t,x) +\frac{1}{2} a(t,x) : D^{2}\phi(t,x) \, ,
\quad\text{for all }
(t,x) \in [0,T]\times \bR^{d} \, ,
\end{align*}
where we define $a(t,x)=\sigma(t,x) \sigma(t,x)^{\intercal}$, $A:B := \text{trace}(AB^{\intercal})$, ${}^\intercal$ is the transpose of a matrix and $D$, $D^{2}$ denotes the usual multi-dimensional spatial differential operators of order one and two (see \cite{Evans1998}).

It well known by the Feynman-Kac formula that if a unique classical solution $v \in C_{b}^{1,2}$ exists to the following PDE,
\begin{equation*}
\begin{cases}
\partial_{t} v(t,x) + \cL v(t,x) =0 \, ,
\\
v(T,x)=g(x) \, ,
\end{cases}
\end{equation*}
for $g$ a Lipschitz continuous function, then the solution of this PDE admits a stochastic representation, $v(t,x)=\bE[g(X_{T})|X_{t}=x]$. Further, by the use of branching diffusions (see \cite{LabordereEtAl2016}) or BSDEs (see \cite{CrisanManolarakis2010}), one is able to obtain a stochastic representation for semi-linear PDEs of the form
\begin{equation*}
\begin{cases}
\partial_{t} v(t,x) + \cL v(t,x) =f(t,x,v,Dv) \, ,
\\
v(T,x)=g(x) \, ,
\end{cases}
\end{equation*}
for $f$ and $g$ nice enough.

\subsection{Heuristic derivation of the idea of our work}
\label{Sec:Heuristics}
Much of the current literature on branching diffusions and regime switching is technical and complex, to aid the presentation of this paper we give an introductory outline of our work.
The ultimate goal of our paper is to construct a stochastic representation of PDEs with only first order spatial derivatives and develop a way to deal with the corresponding 2nd order nonlinearity. We consider PDEs of the form
\begin{equation}
\label{Eq:Transport PDE}
\begin{cases}
\partial_{t} v(t,x) + b(t,x) \cdot D v(t,x) =0 \, ,
\\
v(T,x)=g(x) \, ,
\end{cases}
\end{equation}
for notational convenience we will work in one spatial dimension here (hence $D= \partial_{x}$). 
The problem with constructing a stochastic representation involving the use of It\^o's formula is that we automatically obtain a second order derivative. However, it is known that arguments from branching diffusion can be used to deal with higher order derivatives through the Bismut-Elworthy-Li formula (automatic differentiation as developed in \cite{FournieEtAl1999}). Let us assume that $v$ solving \eqref{Eq:Transport PDE} is a unique classical solution which is $C_{b}^{1,2}$ (i.e. we can apply It\^{o}'s formula to $v$), then we can consider the following equivalent PDE
\begin{equation*}
\begin{cases}
\partial_{t} v(t,x) + b(t,x)\partial_{x} v(t,x) 
+ \frac{1}{2}\sigma_{0}^{2} \partial_{xx} v(t,x) 
-\frac{1}{2}\sigma_{0}^{2} \partial_{xx} v(t,x) =0 \, ,
\\
v(T,x)=g(x) \, ,
\end{cases}
\end{equation*}
where $\sigma_{0}$ is some constant. In fact, as considered in \cite{LabordereEtAl2015}, we can consider the equivalent PDE,
\begin{equation}
\label{Eq:Transport PDE Altered}
\begin{cases}
\partial_{t} v(t,x) + b_{0}\partial_{x} v(t,x) 
                    + \frac{1}{2}\sigma_{0}^{2} \partial_{xx} v(t,x) 
 + \big( b(t,x) - b_{0} \big) \partial_{x} v(t,x) -\frac{1}{2}\sigma_{0}^{2} \partial_{xx} v(t,x) =0 \, ,
\\
v(T,x)=g(x) \, ,
\end{cases}
\end{equation}
where $b_{0}$ is also some constant. 

\textbf{Stochastic Representation.}
Using the Feynman-Kac formula one can easily obtain the following stochastic representation of the solution to \eqref{Eq:Transport PDE Altered},
\begin{equation}
\label{Eq:Basic Feynman-Kac}
v(t,x)= \bE \left[ g(\bar{X}_{T}) 
+ \int_{t}^{T} 
\Big( \big( b(s, \bar{X}_{s}) - b_{0} \big) \partial_{x} v(s, \bar{X}_{s})  
		      -\frac{1}{2}\sigma_{0}^{2} \partial_{xx}v(s, \bar{X}_{s})
		      \Big) \dd s
~ \Big | 
~ \bar{X}_{t}=x \right] \, ,
\end{equation}
where the driving SDE satisfies
\begin{equation}
\label{Eq:Driving SDE}
\dd \bar{X}_{s}= b_{0} \dd s + \sigma_{0} \dd W_{s} \, , \quad \bar{X}_{t}=x \,\quad s\in[t,T] .
\end{equation}
One can observe that such a representation holds provided our constants are $\cF_{t}$ measurable. 

\textbf{Introduce a new random variable.}
Following a standard branching diffusion style argument, alongside the Brownian motion, $W$, we also consider an independent random variable $\tau$ with density $f>0$ on $[0, T-t+ \epsilon]$ for $\epsilon >0$ and denote by $\overline{F}$ the corresponding survival function, namely for $s \in \bR_{+}$ $\overline{F}(s) := \int_{s}^{\infty} f(r) \dd r$. Consider some nice functions $\psi$ and $\phi$, then following representation holds
\begin{align*}
\psi(T) + \int_{t}^{T} \phi(s) \dd s 
& = 
\frac{\psi(T) \overline{F}(T-t)}{\overline{F}(T-t)} + \int_{t}^{T} \frac{1}{f(s-t)} \phi(s) f(s-t) \dd s
\\
& =
\bE_{f}\left[\1_{\{\tau \ge T-t\}}\frac{\psi(T)}{\overline{F}(T-t)} + \1_{\{\tau < T-t\}}\frac{1}{f(\tau)} \phi(t+\tau) \right],
\end{align*}
where $\bE_{f}$ denotes the expectation for the random variable $\tau$.

\textbf{Rewriting the stochastic representation \eqref{Eq:Basic Feynman-Kac}.}
Applying this to the Feynman-Kac representation \eqref{Eq:Basic Feynman-Kac} yields,
\begin{align}
\label{Eq:Feynman-Kac PDE}
v(t,x)
= & \bE \Bigg[ \frac{g(\bar{X}_{T})}{\overline{F}(T-t)}\1_{\{t+\tau \ge T\}} 
+~ \1_{\{t+\tau < T\}} \frac{1}{f(\tau)}
\Big[ -\frac{1}{2}\sigma_{0}^{2} \partial_{xx}v(t+\tau, \bar{X}_{t+\tau})
 \notag
\\ 
&
\hspace{3cm}
+ \big( b(t+\tau, \bar{X}_{t+\tau}) - b_{0} \big) \partial_{x} v(t+\tau, \bar{X}_{t+\tau})  
\Big]
~ \Big | 
~ \bar{X}_{t}=x \Bigg] \, .
\end{align}
One may note the abuse of notation here, the original Feynman-Kac representation expectation was only w.r.t. the Brownian motion, while \eqref{Eq:Feynman-Kac PDE} is w.r.t. both $\tau$ and the Brownian motion.
To make the notation easier we now introduce the following stochastic sequence of times (stochastic mesh on the interval $[t,T]$), $t=:T_{0}<T_{1}< \dots < T_{N_{T}}< T_{N_{T}+1}:=T$ constructed as follows, take a sequence of i.i.d. copies of $\tau$, then set $T_{k+1}=(T_{k}+\tau^{(k)}) \wedge T$ for $k \in \Lambda \subset \bN$, where $\Lambda$ is the set of integers (of stochastic length) $\{1, \dots, N_{T}+1\}$. Using this mesh we then define $\Delta T_{k+1}=T_{k+1}-T_{k}$ and $\Delta W_{T_{k+1}}= W_{T_{k+1}}-W_{T_{k}}$. 

\textbf{Choosing the SDE's coefficients.} 
Let us now consider a good choice of constant for $b_{0}$ (we define $\sigma_{0}$ later). As discussed in \cites{LabordereEtAl2015,DoumbiaEtAl2017}, one can use the so called \emph{frozen coefficient} function which defines the Euler scheme. That is, we may define the SDE $\bar{X}$ recursively over the random mesh by 
\begin{equation}
\label{Eq:Euler Recursion}
\bar{X}_{T_{k}}=
\bar{X}_{T_{k-1}} + b(T_{k-1}, \bar{X}_{T_{k-1}})\Delta T_{k} + \sigma_{k-1} \Delta W_{T_{k}}, \quad \bar{X}_{0}=x \, ,
\end{equation}
for $k \in \Lambda$.
Define $\theta_{k-1}$ as the times in the mesh and position of the SDE up to time $T_{k-1}$ i.e. $\theta_{k-1}:= (T_{1}, \dots, T_{k-1}, x, \bar{X}_{T_{1}}, \dots, \bar{X}_{T_{k-1}})$. Furthermore define the functions $\bar{b}(\theta_{k-1},s,\bar{X}_{s})=b(T_{k-1}, \bar{X}_{T_{k-1}})$ and $\sigma(\theta_{k-1},s)=\sigma_{k-1}$ for $T_{k-1} < s $. Then the SDE defined recursively by,
\begin{equation}
\label{Eq:Euler Integral Recursion}
\bar{X}_{T_{k}}=
\bar{X}_{T_{k-1}} + \int_{T_{k-1}}^{T_{k}}\bar{b}(\theta_{k-1},s,\bar{X}_{s}) \dd s + \int_{T_{k-1}}^{T_{k}}\sigma(\theta_{k-1},s) \dd W_{s} \, ,
\end{equation}
 is the Euler scheme in \eqref{Eq:Euler Recursion}. Moreover, it is clear that the coefficients $\bar{b}(\theta_{k}, \cdot)$ and $\sigma(\theta_{k}, \cdot)$ are $\cF_{T_{k}}$-adapted, hence can be used in \eqref{Eq:Feynman-Kac PDE}. Using the coefficients coming from the Euler scheme is key here since we can simulate an Euler scheme exactly and hence the SDE appearing in \eqref{Eq:Feynman-Kac PDE} can be simulated exactly (which leads to the unbiased representation).

\begin{remark}
	We draw attention to a subtlety in the notation, we will define $\sigma$ on intervals of the form $(\cdot, \cdot]$, thus $\sigma$ is constant over each interval in the time mesh (as is the case in the Euler scheme).
\end{remark}

\textbf{Obtaining a representation for the derivatives.} 
The only terms left to consider in \eqref{Eq:Feynman-Kac PDE} are the derivatives of $v$.
We will formulate rigorous results in Section \ref{Sec:Main Results}, for now let us assume that all functions are sufficiently smooth and with good properties. We construct the Bismut-Elworthy-Li formula (automatic differentiation) w.r.t. the SDE \eqref{Eq:Euler Integral Recursion}.
From \cite{FournieEtAl1999}*{Assumption 3.1} the following integration by parts relation holds for any square integrable function $\phi$,
\begin{align*}
\partial_{x} \bE[\phi(X_{s})|X_{t}=x]
=
\bE \left[ \phi(X_{s}) \int_{t}^{s} \sigma(u)^{-1}Y(u) \mu(u) \dd W_{u} ~ \Big| ~ X_{t}=x
\right] \, ,
\end{align*}
where $Y$ is the first variation process of the SDE $X$ and $\mu$ is any function such that $\int_{t}^{s} \mu(u) \dd u =1$. In the case of the SDE being \eqref{Eq:Euler Integral Recursion}, it is clear that the first variation process is constant equal to one (note $\sigma$ does not have a space dependence). Typically one takes constant $\mu=1/(s-t)$, thus for \eqref{Eq:Euler Integral Recursion} we obtain,
\begin{align*}
\partial_{x} \bE[\phi(X_{T_{1}})|X_{t}=x]
=
\bE \left[ \phi(X_{T_{1}}) \frac{1}{\Delta T_{1}} \int_{t}^{T_{1}} \sigma(\theta_{0}, u)^{-1} \dd W_{u} ~ \Big| ~ X_{t}=x
\right] \, ,
\end{align*}
The same method yields a similar expression for the second derivative
\begin{align*}
\partial_{xx} \bE[\phi(X_{T_{1}})|X_{t}=x]
=
\bE \left[  \frac{\phi(X_{T_{1}})}{\Delta T_{1}^{2}} 
\left(
\left(\int_{t}^{T_{1}} \sigma(\theta_{0}, u)^{-1} \dd W_{u} \right)^2
-
\int_{t}^{T_{1}} (\sigma(\theta_{0}, u)^{-1})^{2} \dd u
\right)
~ \Big| ~ X_{t}=x
\right] \, ,
\end{align*}
From this result and using the fact that $\sigma$ is constant between mesh points we obtain for the second derivative
\begin{align*}
\partial_{xx} v(t,x) 
=& 
\bE \Bigg[   \frac{g(X_{T})}{\overline{F}(\Delta T_{1})}\1_{\{ T_{1} \ge T\}} \frac{1}{\Delta T_{1}^{2}} 
\left(
\left(\sigma(\theta_{0}, T_{1})^{-1} \Delta W_{T_{1}} \right)^2
-
(\sigma(\theta_{0},T_{1})^{-1})^{2} \Delta T_{1}
\right)
\\ &
+  \frac{ \1_{\{T_{1} < T\}} }{f(\Delta T_{1})}
\left(
\big( b(T_{1}, \bar{X}_{T_{1}}) - \bar{b}(\theta_{0},T_{1},\bar{X}_{T_{1}}) \big) 
\partial_{x} v(T_{2}, \bar{X}_{T_{2}})  
-
\frac{1}{2}\sigma(\theta_{0}, T_{1})^{2} \partial_{xx}v(T_{1}, X_{T_{1}})  
\right)
\\
&
\quad
\times
\frac{\left(\sigma(\theta_{0}, T_{1})^{-1} \Delta W_{T_{1}} \right)^2
	-
	(\sigma(\theta_{0}, T_{1})^{-1})^{2} \Delta T_{1}}
{\Delta T_{1}^{2}}
~ \Big | 
~ X_{t}=x \Bigg] \, ,
\end{align*}
the $\partial_{x}v$ term is similar. The idea of branching diffusion style algorithms is to continuously substitute in terms involving the solution until we remove the dependence on it. 
Of course, $v(t,x)$ does not appear inside the expectation, however, by using the tower property and flow property of the SDE we are able to derive the corresponding representations for $\partial_{x} v(T_{k},X_{T_{k}})$ and $\partial_{xx} v(T_{k},X_{T_{k}})$.

\textbf{Rewriting the stochastic representation \eqref{Eq:Feynman-Kac PDE}.}
Substituting in the expressions for $\partial_{x} v(T_{1},X_{T_{1}})$ and $\partial_{xx} v(T_{1},X_{T_{1}})$ into \eqref{Eq:Feynman-Kac PDE} yields,
\begin{align*}
& v(t,x)
\\
& = \bE \Bigg[ \frac{g(\bar{X}_{T})}{\bar{F}(\Delta T_{1})}\1_{\{T_{1} \ge T\}} 
+ 
\1_{\{T_{1} < T\}} \frac{1}{f(\Delta T_{1})}
\bE \Bigg[ \overline{\cW}_{2}
\Bigg\{
\frac{g(\bar{X}_{T})}{\overline{F}(\Delta T_{2})}\1_{\{T_{2} \ge T\}}
+
\1_{\{T_{2} < T\}} \frac{1}{f(\Delta T_{2})}
\\
&
\times \left[ 
\big( b(T_{2}, \bar{X}_{T_{2}}) - \bar{b}(\theta_{1},T_{2},\bar{X}_{T_{2}}) \big) 
\partial_{x} v(T_{2}, \bar{X}_{T_{2}})  
-
\frac{1}{2}\sigma(\theta_{1},T_{2})^{2} \partial_{xx}v(T_{2}, \bar{X}_{T_{2}})
\right]
\Bigg\}
~ \Bigg | ~ \bar{X}_{T_{1}}
\Bigg]
~ \Bigg | 
~ \bar{X}_{t}=x \Bigg] \, ,
\end{align*}
where $\overline{\cW}_{k}$ is the so-called Malliavin weight stemming from the automatic differentiation,
\begin{align*}
\overline{\cW}_{k}:=
\frac{ b(T_{k-1}, \bar{X}_{T_{k-1}}) - \bar{b}(\theta_{k-2},T_{k-1},\bar{X}_{T_{k-1}}) }
{\sigma(\theta_{k-1},T_{k})}
\frac{\Delta W_{T_{k}}}{ \Delta T_{k}}
-
\frac{1}{2} \frac{\sigma(\theta_{k-2},T_{k-1})^{2}}{\sigma(\theta_{k-1},T_{k})^{2} } \left(
\frac{\Delta W_{T_{k}}^{2} - \Delta T_{k}}{\Delta T_{k}^{2}}
\right). 
\end{align*}
One observes that this Feynman-Kac representation now only depends on the solution $v$ if $T_{2}<T$. 

\textbf{Taking the limit.}
Following the standard procedure in branching diffusions (see \cites{Labordere2012,LabordereEtAl2014,LabordereEtAl2016}), executing the same argument multiple times removes the dependence on $v$ on the right hand side. Following \cite{DoumbiaEtAl2017} we introduce the following notation,
\begin{align}
\label{Eq:Explicit Malliavin Weights}
M_{k+1}=\Delta b_{k} \sigma(\theta_{k},T_{k+1})^{-1} \frac{\Delta W_{T_{k+1}}}{ \Delta T_{k+1}} 
\quad 
\text{ and }
\quad
V_{k+1}= -\frac{1}{2} \frac{\sigma(\theta_{k-1},T_{k})^{2}}{\sigma(\theta_{k},T_{k+1})^{2} } 
\left(
\frac{\Delta W_{T_{k+1}}^{2} - \Delta T_{k+1}}{\Delta T_{k+1}^{2}}
\right) \, .
\end{align}
where $\Delta b_{k}= b(T_{k},\bar{X}_{T_{k}})-\bar{b}(\theta_{k-1},T_{k}, \bar{X}_{T_{k}})= b(T_{k}, \bar{X}_{T_{k}})-b(T_{k-1}, \bar{X}_{T_{k-1}})$. Further define the terms 
\begin{align}
\label{Eq:Malliaivn Weight Term}
P_{k+1}:= \frac{M_{k+1}+\frac{1}{2}V_{k+1}}{f(\Delta T_{k})} \quad \text{for} ~~ k \in \Lambda.
\end{align}
It is then clear that the solution to the PDE can be written as follows,
\begin{align}
\label{Eq:Non L2 representation}
v(t,x)=\bE \left[ \frac{g(\bar{X}_{T_{N_{T}+1}})}{\bar{F}(\Delta T_{N_{T}+1})}
\prod_{k=2}^{N_{T}+1} P_{k}
~
\Bigg | ~
\bar{X}_{t}=x
\right] \, .
\end{align}
Although this relation is useful for us, in its current form it is not square integrable, thus we need to use some variance reduction techniques in order to use Monte Carlo. Moreover, many of the operations above require some form of integrability, these points will be the main focus of the next section.



\section{Stochastic Representation for a toy transport PDE}
\label{Sec:Main Results}

The goal of the paper is to derive a square-integrable representation that solves a PDE of the form, 
\begin{equation}
	\label{Eq:High Dimension Transport PDE}
	\begin{cases}
		\partial_{t} v(t,x) + b(t)\cdot Dv(t,x) =0  \quad\text{for all }
		(t,x) \in [0,T)\times \bR^{d} \, ,
		\\
		v(T,x)=g(x) \, .
	\end{cases}
\end{equation}
\begin{remark}
	We show the representation in the case $b:[0,T]\to \bR$ is independent of space. This ensures finite variance, we shall return to the case of space dependency later.
\end{remark}

We wish to consider SDEs of the form \eqref{Eq:Euler Integral Recursion}, in $d$-dimensions this is,
\begin{align*}
	\dd \bar{X}_{s}= \bar{b}(\theta,s) \dd s + \sigma(\theta,s) \bI_{d} \dd W_{s} \, , \quad \text{for } s \in [t,T] ~ \text{and } \bar{X}_{t}=x \, , 
\end{align*}
where $\bI_{d}$ is the $d$-dimensional identity matrix.
Unlike typical stochastic representations, $\sigma$ is not fixed by the PDE, thus we have the freedom to choose $\sigma$. Although, the representation is somewhat independent of the precise choice of $\sigma$, the variance of the estimate (and hence the usefulness) heavily depends on $\sigma$.

In order to keep our representation and in particular our proofs as readable as possible, we consider only the one dimensional case. As one can clearly see though, due to fact that $\sigma$ is a scalar multiplied by the identity, all our arguments generalise to the higher dimensional case. Of course as $\sigma$ is not fixed in this case, it may be that other representations especially in high dimension may yield superior results. However, our goal here is purely to obtain a representation with finite variance.

The previous section outlined how one builds the stochastic representation without going into detail about when the various steps are applicable. We now want to show that this representation holds under some integrability and regularity assumptions. 
In the previous section we required two types of random variable, namely a driving Brownian motion and an i.i.d. sequence of random times $\tau^{(k)}$ with density $f$, independent of the Brownian motion and $k \in \Lambda$ as before. Thus consider the probability space $(\Omega, \cF, \bP)$ generated by these random variables, we also denote by $\bP_{W}$ and $\bP_{f}$ the probability measure ($\bE_{W}$ and $\bE_{f}$ the corresponding expectation) restricted to the Brownian motion and random times respectively. With this notation, one may think of $\bP$ as the product measure $\bP_{W} \otimes \bP_{f}$. 
The corresponding filtration $\cF_{t}$ is the sigma-algebra generated by the set of random times up to $t$ i.e. $\max \{k: T_{k} \le t\}$ and the Brownian motion up to $t$, hence, $\cF_{t}:= \sigma(T_{1}, \dots, T_{k}, (W_{s})_{s \le t})$.

Let us first state the assumptions we will use.
\begin{assumption}
	\label{Assume:Main SDE Assumptions}
	We assume the drift, $b$ is uniformly Lipschitz in time.
\end{assumption}
The analysis we carry out using regime switching techniques is sufficiently difficult to present that we assume the existence of a good enough solution to the transport PDE, as opposed to assuming sufficient conditions that would allow us to derive the said solution. Waiving the next assumption is left for future work.
\begin{assumption}
	\label{Assume:Bounded g}
	Firstly we assume that there exists a unique solution $v \in C_{b}^{1,3}([0,T], \bR^{d})$ to \eqref{Eq:High Dimension Transport PDE}. In particular, we have that the terminal condition function $g$ of the PDE satisfies $g \in C_{b}^{2}$.
\end{assumption}
The assumption on $g$ is not necessary since it follows from $v \in C_{b}^{1,3}$, however, we make this explicit since it is all we require for our estimator to be of finite variance.
It is possible to put some conditions on $b$ and $g$ leading to a unique solution for general transport PDEs see \cite{Kato1975} for example. We do not go into detail here as this will again be the subject of future work.

We consider the particles to have a life time given by Gamma distributed random variables, i.e. $\tau$ has density,
\begin{align}
\label{Eq:Gamma Density}
f(s) := f_{\Gamma}^{\kappa, \eta}(s) = \frac{s^{\kappa-1} \exp(-s/\eta)}{\Gamma(\kappa) \eta^{\kappa}} \, , \quad \text{for all }~ s >0 
~ ~ \text{where } ~ \kappa, ~ \eta >0 \, ,
\end{align}
where $\Gamma$ is the Euler function $\Gamma(y)=\int_{0}^{\infty}x^{y-1} \exp(-x) \dd x$.

 We will use a \emph{mesh dependent} coefficient for $\sigma$ relying on the times at which the regime switching occurs,
\begin{align}
\label{Eq:Sigma Defn}
\sigma(\theta_{k-1},s):= \sigma_{0} \prod_{i=1}^{k-1} \Delta T_{i}^{n} \quad \text{for } s \in (T_{k-1},T_{k}] ~ ,  ~
k=1,\dots, N_{T}+1 \, ,
~~ 
n \in \bR
~ ~ \text{and}
~ ~ \sigma_{0} \in \bR_{+} \, ,
\end{align} 
hence $\sigma(\theta_{k-1}, T_{k})= \sigma_{0} \prod_{i=1}^{k-1} \Delta T_{i}^{n}$, with the convention $\prod_{i=1}^{0} \cdot =1$.

\begin{remark}[Adaptedness of $\sigma$]
	Even though our $\sigma$ depends on the stochastic mesh, it is $\cF_t$-adapted. This is of fundamental importance to show that the estimator in \eqref{Eq:Variance Reduced Representation} solves the PDE \eqref{Eq:High Dimension Transport PDE}.
\end{remark}

 We make an assumption on the parameters of $\sigma$ and $f$.
\begin{assumption}
	\label{Assume:Finite Variance}
	The power exponent $n$ in the diffusion coefficient \eqref{Eq:Sigma Defn} satisfies $n \le -1$. The shape parameter of the Gamma random variable, \eqref{Eq:Gamma Density}, is $\kappa =1/2$.
\end{assumption}

\begin{remark}
	Under Assumption \ref{Assume:Finite Variance}, $\sigma$ is a positive function bounded from below away from zero\footnote{To see this, note that $n<0$, hence for $\sigma$ to be zero, we require a set of $\Delta T_{k} \ge 0$ for $k=1, \dots, N_{T+1}$, such that $\sum_{k=1}^{N_{T+1}} \Delta T_{k} =T$ and $\prod_{k=1}^{N_{T+1}} \Delta T_{k} = \infty$. Which clearly does not exist.}. The bounds on $n$ and $\kappa$ are mainly for convenience in order for the proof of Proposition \ref{Prop:Representation} to follow.
\end{remark}
As was alluded to in Section \ref{Sec:Representation}, \eqref{Eq:Non L2 representation} was not useful since it did not have finite second moment. To solve this problem we employ variance reduction techniques, namely antithetic variables and control variates. Consider the following auxiliary random variables, $\beta:= (\beta_{1} + \beta_{2})/2$ with
\begin{equation}
\label{Eq:Beta Terms}
\begin{cases}
\beta_{1}:= 
\dfrac{g(\bar{X}_{T_{N_{T}+1}})-g(\bar{X}_{T_{N_{T}}} + b(T_{N_{T}}) \Delta T_{N_{T}+1}) }{\overline{F}(\Delta T_{N_{T}+1})}
\dfrac{M_{N_{T}+1}+\frac{1}{2}V_{N_{T}+1}}{f(\Delta T_{N_{T}})} \, ,
\\
\beta_{2}:= 
\dfrac{g(\hat{X}_{T_{N_{T}+1}})-g(\bar{X}_{T_{N_{T}}} + b(T_{N_{T}}) \Delta T_{N_{T}+1}) }{\overline{F}(\Delta T_{N_{T}+1})}
\dfrac{-M_{N_{T}+1}+\frac{1}{2}V_{N_{T}+1}}{f(\Delta T_{N_{T}})}
\, ,
\end{cases}
\end{equation}
where $\hat{X}$ is the antithetic random variable associated to $\bar{X}$ i.e. the Euler scheme defined by, $\hat{X}_{T_{k}}=\bar{X}_{T_{k-1}} + b(T_{k-1})\Delta T_{k} - \sigma(\theta_{k-1}, T_{k}) \Delta W_{T_{k}}$ and $V$ and $M$ as defined in \eqref{Eq:Explicit Malliavin Weights}. It is straightforward to see that the additional $g$ term is a control variate since its input is independent of Brownian motion $\Delta W_{T_{N_{T}+1}}$.
One can further understand $(\beta_{1},\beta_{2})$ as an antithetic pair.


We now state our main result of the paper.

\begin{theorem}
	\label{Thm:Representation Solves the PDE}
	[Representation Solves the PDE]
	Let Assumptions \ref{Assume:Main SDE Assumptions}, \ref{Assume:Bounded g} and \ref{Assume:Finite Variance} hold, and let us denote by $\hat{v}: [0,T] \times \bR \rightarrow \bR$ the following function,
	\begin{align}
	\label{Eq:Variance Reduced Representation}
	\hat{v}(t,x):= \bE \left[\beta \prod_{k=2}^{N_{T}} P_{k} \1_{\{N_{T} \ge 1\}}
	~
	\Bigg | ~
	\sigma(\theta_{0},t), ~
	X_{t}=x \right] 
	+
	\bE \left[ \frac{g(\bar{X}_{T_{1}})}{\overline{F}(\Delta T_{1})} \1_{\{N_{T}=0\}}
	~
	\Bigg | ~
	\sigma(\theta_{0},t), ~
	X_{t}=x \right], 
	\end{align}
	with $\{P_{k}\}_{k}$ as defined in \eqref{Eq:Malliaivn Weight Term}.
	Then $\hat{v}$ solves the PDE \eqref{Eq:High Dimension Transport PDE}, namely $\hat{v}=v$ (hence $\hat{v}$ is an unbiased estimator of $v$).
	Moreover, the stochastic process generating $\hat{v}$ is square integrable and hence of finite variance.
\end{theorem}

\paragraph{Outline of proof}
The proof of Theorem \ref{Thm:Representation Solves the PDE} requires several steps which we show in the following order.
\begin{enumerate}
	\item Take $\tilde{v}$ in \eqref{Eq:Variance Reduced Representation}, which is the expected value of a stochastic process (estimator).
	\item Show that the estimator is square integrable, Proposition \ref{Prop:Representation}.
	\item Show that under enough integrability a stochastic representation to \eqref{Eq:High Dimension Transport PDE} exists when a solution in $C_{b}^{1,3}([0,T], \bR)$  exists, Theorem 	\ref{theo:StochasticRepresentation}.
	\item Show that \eqref{Eq:Variance Reduced Representation}, satisfies the integrability conditions in Theorem \ref{theo:StochasticRepresentation} and thus solves \eqref{Eq:High Dimension Transport PDE}, Theorem \ref{Thm:Integrability Hold}.
\end{enumerate}



\subsection[Variance analysis for a specific diffusion coefficient]{Variance analysis for a specific diffusion coefficient}
 Since our regime switching algorithm does not create new particles, our computational complexity for any Monte Carlo realisation is only $O(C (N_{T}+1))$, since $T < \infty$, it is clear we have finite computational complexity. We therefore only need to consider the variance of the estimator.
 We obtain the following.
\begin{proposition}
	\label{Prop:Representation}
	Let Assumptions \ref{Assume:Main SDE Assumptions}, \ref{Assume:Bounded g} and \ref{Assume:Finite Variance} hold. Then the random variable appearing in \eqref{Eq:Variance Reduced Representation},
	\begin{align*}
	\beta \prod_{k=2}^{N_{T}} P_{k} \1_{\{N_{T} \ge 1\}}
	+
	\frac{g(\bar{X}_{T_{1}})}{\overline{F}(\Delta T_{1})} \1_{\{N_{T}=0\}}
	\end{align*}	
has finite variance. 
\end{proposition} 

Although this proof is argued in a similar style to the proof of Proposition 4.1 in \cite{DoumbiaEtAl2017}, there are many subtle differences and we overall require a more refined analysis of the various terms to ensure our estimator has finite second moment. We point in particular to the ``Interval splitting'' argument in order to deal with instability in the last time point of the random mesh. This is essential to deal with the second order term that appears.
\begin{proof}[Proof. {[Finite variance of the estimator]}]
	Consider $\overline{\cF}_{k}$ the sigma-algebra generated by the set of random times up to $T_{k+1}$ and the Brownian motion up to $T_{k}$, hence\footnote{One should note the small but critical distinction between $\cF_{t}$ and $\overline{\cF}_{k}$.}, $\overline{\cF}_{k}:= \sigma(T_{1}, \dots, T_{k+1}, (W_{s})_{s \le T \wedge T_{k}})$.
	
	Throughout the proof, for ease of writing we suppress the condition in the expectation of the process starting at $x$ at time $t$.
	
	In order to show finite variance we only need to show finite second moment (the dominant term), further note that due to the indicators we obtain no cross term.
	Looking first at the second term of \eqref{Eq:Variance Reduced Representation}, by the bounds on the coefficients on the SDE and the Lipschitz property of $g$ we have $\bE[g(\bar{X}_{T_{1}})^{2}]<\infty$, and $\overline{F}(T-t)>0$, thus we have finite variance on the second term. For the first term in \eqref{Eq:Variance Reduced Representation}, we can rewrite the second moment as,
	\begin{align*}
	\bE \left[ \left( \beta \prod_{k=2}^{N_{T}} P_{k} \right)^{2} \1_{\{N_{T} \ge 1\}} \right] 
	=
	\sum_{\ell =1}^{\infty} \bE \left[ \left( \beta \prod_{k=2}^{N_{T}} P_{k} \right)^{2} \Bigg | N_{T}=\ell \right] \times \bP[N_{T}=\ell] \, .
	\end{align*}
	
	In order to tackle this term we split the proof into several steps by bounding various quantities then combining them together to show the sum is bounded.
	We also note that we often work with conditional expectations, hence statements involving them are to be understood in the $\bP$-a.s. sense.

	
	
	
	

	\emph{Step 1: Bounding $\bE[\beta^{2}| \overline{\cF}_{N_{T}}, N_{T}= \ell]$, for $\beta$ from \eqref{Eq:Beta Terms}}.
	\label{Step: Beta}As is standard practice when we only care about showing an estimate to be finite we use $C$ to denote some finite constant which can change over inequalities but crucially can only depend on ``known'' constants such as $T$ etc. By the tower property we can rewrite any term in the sum as,
	\begin{align*}
	\bE \left[ \left( \beta \prod_{k=2}^{N_{T}} P_{k} \right)^{2} \Bigg | N_{T}=\ell \right] 
	=
	\bE \left[ \bE \big[ \beta^{2} | \overline{\cF}_{N_{T}}, N_{T}=\ell \big] \prod_{k=2}^{N_{T}} P_{k}^{2} \Bigg | N_{T}=\ell \right] 
	\, .
	\end{align*}
	Rewriting $\beta$ with $M_{N_{T}+1}$ and $V_{N_{T}+1}$ as common factors then using Young's inequality we obtain,
	\begin{align*}
	& \bE\big[\, \beta^{2} | \bar{\cF}_{N_{T}}, N_{T}=\ell\, \big]
	\\
	&
	\le  
	C \bE \left[ 
	\dfrac{ \left( g(\bar{X}_{T_{N_{T}+1}})- g(\hat{X}_{T_{N_{T}+1}}) \right)^{2}
	}{
	\overline{F}(\Delta T_{N_{T}+1})^{2}
}
\dfrac{M_{N_{T}+1}^{2}}{f(\Delta T_{N_{T}})^{2}}
\Bigg| \overline{\cF}_{N_{T}}, N_{T}=\ell \right]
\\ &
+
C \bE \left[ 
\dfrac{ \left( g(\bar{X}_{T_{N_{T}+1}})+ g(\hat{X}_{N_{T}+1}) -2g(\bar{X}_{T_{N_{T}}} + b(T_{N_{T}}) \Delta T_{N_{T}+1}) \right)^{2}
}{
\overline{F}(\Delta T_{N_{T}+1})^{2}
}
\dfrac{\frac{1}{2}V_{N_{T}+1}^{2}}{f(\Delta T_{N_{T}})^{2}}
\Bigg| \overline{\cF}_{N_{T}}, N_{T}=\ell \right] \, .
\end{align*}
Considering the first term on the RHS, we note by the Lipschitz property of $g$ that,
\begin{align*}
|g(\bar{X}_{T_{N_{T}+1}})-g(\hat{X}_{T_{N_{T}+1}})| 
\le
L| \bar{X}_{T_{N_{T}+1}}-\hat{X}_{T_{N_{T}+1}} |
\le
C |  \sigma(\theta_{N_{T}},T_{N_{T}+1}) \Delta W_{T_{N_{T}+1}}| \, .
\end{align*}
Hence using this bound and the representation for $M_{N_{T}+1}$ (see \eqref{Eq:Explicit Malliavin Weights}),
\begin{align*}
\bE & \left[ 
\dfrac{ \left( g(\bar{X}_{T_{N_{T}+1}})- g(\hat{X}_{T_{N_{T}+1}}) \right)^{2}
}{
\overline{F}(\Delta T_{N_{T}+1})^{2}
}
\dfrac{M_{N_{T}+1}^{2}}{f(\Delta T_{N_{T}})^{2}}
\Bigg| \overline{\cF}_{N_{T}}, N_{T}=\ell \right]
\\
& \le
C 
\dfrac{
	\Delta b_{N_{T}}^{2}
}{
f(\Delta T_{N_{T}})^{2}}
\bE \left[ 
\left(\Delta W_{T_{N_{T}+1}}\sigma(\theta_{N_{T}},T_{N_{T}+1})
\right)^{2}
\left( \frac{\Delta W_{T_{N_{T}+1}}}{\Delta T_{N_{T}+1}} \sigma(\theta_{N_{T}},T_{N_{T}+1})^{-1} \right)^{2}
\Bigg| \overline{\cF}_{N_{T}}, N_{T}=\ell \right] 
\\
& =
C 
\dfrac{
	\Delta b_{N_{T}}^{2}
}{
f(\Delta T_{N_{T}})^{2}}
\, ,
\end{align*}
where we used $1/\overline{F}(\Delta T_{N_{T}+1})^{2} \le C$ in the inequality.
For the second term on the RHS, it is more complex, let us first split the terms using Cauchy-Schwarz,
\begin{align*}
&\bE \left[ 
\left( g(\bar{X}_{T_{N_{T}+1}})+ g(\hat{X}_{N_{T}+1}) -2g(\bar{X}_{T_{N_{T}}} + b(T_{N_{T}}) \Delta T_{N_{T}+1}) \right)^{2}
V_{N_{T}+1}^{2}
\Bigg| \overline{\cF}_{N_{T}}, N_{T}=\ell \right] 
\\
&
\le
\bE \left[ 
\left( g(\bar{X}_{T_{N_{T}+1}})+ g(\hat{X}_{N_{T}+1}) -2g(\bar{X}_{T_{N_{T}}} + b(T_{N_{T}}) \Delta T_{N_{T}+1}) \right)^{4}
\Bigg| \overline{\cF}_{N_{T}}, N_{T}=\ell 
\right]^{1/2}
\\
& \qquad
\times
\bE \left[ 
V_{N_{T}+1}^{4}
\Bigg| \overline{\cF}_{N_{T}}, N_{T}=\ell \right]^{1/2}
\, .
\end{align*}
Let us firstly focus on the $g$ term. Consider the ODE on the interval $s \in [T_{N_{T}}, T_{N_{T}+1}]$,
\begin{align*}
\frac{\dd Y_{s}}{\dd s} =  b(T_{N_{T}}) \, ,
\qquad
Y_{T_{N_{T}}}= \bar{X}_{T_{N_{T}}} \, .
\end{align*}
Then, the solution is  $Y_{T_{N_{T}+1}} = \bar{X}_{T_{N_{T}}} + b(T_{N_{T}}) \Delta T_{N_{T}+1}$. Consequently,
\begin{align}
g \bigg( & \bar{X}_{T_{N_{T}}} + b(T_{N_{T}}) \Delta T_{N_{T}+1} \bigg)
-
g(\bar{X}_{T_{N_{T}}})
\notag
\\
& =
\int_{T_{N_{T}}}^{T_{N_{T}+1}} g'(Y_{s}) \dd Y_{s} 
=
\int_{T_{N_{T}}}^{T_{N_{T}+1}} g' \Big( \bar{X}_{T_{N_{T}}} + b(T_{N_{T}})(s-T_{N_{T}}) \Big) b\big(T_{N_{T}} \big) \dd s \, .
\label{Eq:Deterministic difference in g}
\end{align}
By applying It\^{o}'s formula to $g(\bar{X}_{T_{N_{T}+1}})$ and $g(\hat{X}_{T_{N_{T}+1}})$ (recall $g \in C_{b}^{2}$), and using \eqref{Eq:Deterministic difference in g} we obtain,
\begin{align}
& g(\bar{X}_{T_{N_{T}+1}}) + g(\hat{X}_{T_{N_{T}+1}}) -2g(\bar{X}_{T_{N_{T}}} + b(T_{N_{T}}) \Delta T_{N_{T}+1}) \notag
\\
& = 
\frac{1}{2} \sigma(\theta_{N_{T}},T_{N_{T}+1})^{2}
\int_{T_{N_{T}}}^{T_{N_{T}+1}} (g''(\bar{X}_{s}) + g''(\hat{X}_{s}))  \dd s 
+
\sigma(\theta_{N_{T}},T_{N_{T}+1})
\int_{T_{N_{T}}}^{T_{N_{T}+1}} (g'(\bar{X}_{s}) - g'(\hat{X}_{s}))  \dd W_{s} \notag
\\
& \quad
+
\int_{T_{N_{T}}}^{T_{N_{T}+1}}  \Big(g'(\bar{X}_{s}) + g'(\hat{X}_{s})- 2g' \Big( \bar{X}_{T_{N_{T}}} + b(T_{N_{T}})(s-T_{N_{T}}) \Big)
\Big)  b(T_{N_{T}}) \dd s \, . 
\label{Eq:Expanded g terms}
\end{align}
Since $g'$ is Lipschitz, we obtain,
\begin{align*}
|g'(\bar{X}_{s}) - g' \Big( \bar{X}_{T_{N_{T}}} + b(T_{N_{T}})(s-T_{N_{T}}) \Big) |
& \le
C | \bar{X}_{s} - \bar{X}_{T_{N_{T}}} + b(T_{N_{T}})(s-T_{N_{T}}) |
\\
& \le
C \sigma(\theta_{N_{T}},T_{N_{T}+1}) |W_{s} - W_{T_{N_{T}}}| \, ,
\end{align*}
the same bound holds for the $g(\hat{X}_{s})$ term.
Thus the following bound can be obtained for the final integral in \eqref{Eq:Expanded g terms}
\begin{align*}
\int_{T_{N_{T}}}^{T_{N_{T}+1}}  & \Big(g'(\bar{X}_{s}) + g'(\hat{X}_{s})- 2g' \Big( \bar{X}_{T_{N_{T}}} + b(T_{N_{T}})(s-T_{N_{T}}) \Big)
\Big)  b(T_{N_{T}}) \dd s 
\\
& \le
C |   b(T_{N_{T}})| \sigma(\theta_{N_{T}},T_{N_{T}+1})
\int_{T_{N_{T}}}^{T_{N_{T}+1}}  |W_{s} - W_{T_{N_{T}}}| \dd s
\, .
\end{align*}
Recalling that we are interested in the fourth moment, using Doob's maximal inequality,
\begin{align*}
\bE & \left[ \left( \int_{T_{N_{T}}}^{T_{N_{T}+1}}  |W_{s} - W_{T_{N_{T}}}| \dd s \right)^{4}
~ \Big| ~ \overline{\cF}_{N_{T}}, N_{T}=\ell
\right]
\\
&
\le
C \Delta T_{N_{T}+1}^{4} \bE \left[ \sup_{T_{N_{T}} \le s \le T_{N_{T}+1}} |W_{s} - W_{T_{N_{T}}}|^{4}
~ \Big| ~ \overline{\cF}_{N_{T}}, N_{T}=\ell
\right]
\le
C \Delta T_{N_{T}+1}^{6} \, .
\end{align*}
For the stochastic integral in \eqref{Eq:Expanded g terms}, again taking the fourth moment we obtain,
\begin{align*}
\sigma(\theta_{N_{T}},T_{N_{T}+1})^{4}
&
\bE \left[ 
\left( 
\int_{T_{N_{T}}}^{T_{N_{T}+1}} (g'(\bar{X}_{s}) - g'(\hat{X}_{s}))  \dd W_{s}
\right)^{4}
~ \Big| ~ \overline{\cF}_{N_{T}}, N_{T}=\ell
\right]
\\
&
=
3
\sigma(\theta_{N_{T}},T_{N_{T}+1})^{4}
\bE \left[  
\left( \int_{T_{N_{T}}}^{T_{N_{T}+1}} (g'(\bar{X}_{s}) - g'(\hat{X}_{s}))^{2}  \dd s \right)^{2}
~ \Big| ~ \overline{\cF}_{N_{T}}, N_{T}=\ell
\right] \, .
\end{align*}
Using that $g'$ is Lipschitz and the difference is given by
\begin{align*}
|g'(\bar{X}_{s}) - g'(\hat{X}_{s})| 
\le 
C |\sigma(\theta_{N_{T}}, T_{N_{T}+1}) (W_{s} - W_{T_{N_{T}}}) 
+
\sigma(\theta_{N_{T}}, T_{N_{T}+1})(W_{s} - W_{T_{N_{T}}})| \, .
\end{align*}
This along with a similar Doob's maximal inequality implies that we can bound the stochastic integral by,
\begin{align*}
\sigma(\theta_{N_{T}},T_{N_{T}+1})^{4}
&
\bE \left[ 
\left( 
\int_{T_{N_{T}}}^{T_{N_{T}+1}} (g'(\bar{X}_{s}) - g'(\hat{X}_{s}))  \dd W_{s}
\right)^{4}
~ \Big| ~ \overline{\cF}_{N_{T}}, N_{T}=\ell
\right]
\le
C \sigma(\theta_{N_{T}},T_{N_{T}+1})^{8} \Delta T_{N_{T}+1}^{4}  \, .
\end{align*}
Recalling that $g''$ is bounded, we can bound the remaining term in \eqref{Eq:Expanded g terms} by a similar term to the stochastic integral to obtain,
\begin{align}
\label{Eq:Ito Bound on g}
\bE \Big[
(g(\bar{X}_{T_{N_{T}+1}}) + g(\hat{X}_{T_{N_{T}+1}}) -2g(\bar{X}_{T_{N_{T}}} + b(T_{N_{T}}) \Delta T_{N_{T}+1})
)^{4} 
& | \overline{\cF}_{N_{T}}, N_{T}=\ell
\Big]
\\
& 
\le
C \sigma(\theta_{N_{T}},T_{N_{T}+1})^{8} \Delta T_{N_{T}+1}^{4} \, . \notag
\end{align}

The above bound was obtained using differentiability and It\^{o}'s formula, however,
it will also be useful for us to note that just using the Lipschitz property yields,
\begin{align*}
\bE \Big[
(g(\bar{X}_{T_{N_{T}+1}}) + g(\hat{X}_{T_{N_{T}+1}}) -2g(\bar{X}_{T_{N_{T}}} + b(T_{N_{T}}) \Delta T_{N_{T}+1})
)^{4} 
& | \overline{\cF}_{N_{T}}, N_{T}=\ell
\Big]
\\
& 
\le
C \sigma(\theta_{N_{T}},T_{N_{T}+1})^{4} \Delta T_{N_{T}+1}^{2} \, .
\end{align*}
Hence we obtain the following stronger bound for the $g$ terms
\begin{align*}
\bE & \Big[
(g(\bar{X}_{T_{N_{T}+1}}) + g(\hat{X}_{T_{N_{T}+1}}) -2g(\bar{X}_{T_{N_{T}}} + b(T_{N_{T}}) \Delta T_{N_{T}+1})
)^{4} 
| \overline{\cF}_{N_{T}}, N_{T}=\ell
\Big]
\\
&
\hspace{5cm}
\le
C \min \left[
\sigma(\theta_{N_{T}},T_{N_{T}+1})^{4} \Delta T_{N_{T}+1}^{2}
,
\sigma(\theta_{N_{T}},T_{N_{T}+1})^{8} \Delta T_{N_{T}+1}^{4}
\right]
\, .
\end{align*}
For the $V$ term,
\begin{align*} 
\bE \left[ V_{N_{T}+1}^{4} | \overline{\cF}_{N_{T}}, N_{T}=\ell \right] 
& \le
C \frac{\sigma(\theta_{N_{T}-1},T_{N_{T}})^{8}}{\sigma(\theta_{N_{T}},T_{N_{T}+1})^{8}}
\frac{1}{\Delta T_{N_{T}+1}^{8}}
\bE \left[ \left(
\Delta W_{T_{N_{T}+1}}^{2}
-
\Delta T_{N_{T}+1}
\right)^{4}
\Big| \overline{\cF}_{N_{T}}, N_{T}=\ell \right] 
\\
&
\le
C \frac{\sigma(\theta_{N_{T}-1},T_{N_{T}})^{8}}{\sigma(\theta_{N_{T}},T_{N_{T}+1})^{8}}
\frac{1}{\Delta T_{N_{T}+1}^{4}}
\, .
\end{align*}
Hence using Cauchy-Schwarz we obtain,
\begin{align*}
& \bE \left[ 
\dfrac{ \left( g(\bar{X}_{T_{N_{T}+1}})+ g(\hat{X}_{T_{N_{T}+1}}) -2g(\bar{X}_{T_{N_{T}}} + b(T_{N_{T}}) \Delta T_{N_{T}+1}) \right)^{2}
}{
\overline{F}(\Delta T_{N_{T}+1})^{2}
}
\dfrac{\frac{1}{2}V_{N_{T}+1}^{2}}{f(\Delta T_{N_{T}})^{2}}
\Bigg| \overline{\cF}_{N_{T}}, N_{T}=\ell \right]
\\
& \le 
C \frac{\sigma(\theta_{N_{T}-1},T_{N_{T}})^{4}}{\sigma(\theta_{N_{T}},T_{N_{T}+1})^{4}}
\frac{1}{\Delta T_{N_{T}+1}^{2}}
\min \left[
\sigma(\theta_{N_{T}},T_{N_{T}+1})^{2} \Delta T_{N_{T}+1}
,
\sigma(\theta_{N_{T}},T_{N_{T}+1})^{4} \Delta T_{N_{T}+1}^{2}
\right]
\dfrac{1}{f(\Delta T_{N_{T}})^{2}}
\, .
\end{align*}
Therefore, the conditional expectation of $\beta^{2}$ can be bounded by,
\begin{align*}
\bE[ \beta^{2}  | \overline{\cF}_{N_{T}}, N_{T}=\ell]
\le  
\dfrac{
	C
}{
f(\Delta T_{N_{T}})^{2}}
\left(
\Delta b_{N_{T}}^{2}
+
\frac{\sigma(\theta_{N_{T}-1},T_{N_{T}})^{4}}{\sigma(\theta_{N_{T}},T_{N_{T}+1})^{2}}
\frac{\min \left[
	1
	,
	\sigma(\theta_{N_{T}},T_{N_{T}+1})^{2} \Delta T_{N_{T}+1}
	\right]}{\Delta T_{N_{T}+1}}
\right)
\, .
\end{align*}






\emph{Step 2: Bounding $ \bE[P_{k+1}^{4}| \overline{\cF}_{k}, N_{T}= \ell]$}.
Let $k \in \Lambda$ and note by Assumption \ref{Assume:Main SDE Assumptions} we obtain,
\begin{align*}
\bE\big[\Delta b_{k}^{4} | \overline{\cF}_{k-1}, N_{T}= \ell\big] 
\le
C \Delta T_{k}^{4} \, .
\end{align*}

From \eqref{Eq:Explicit Malliavin Weights} we observe the following,
\begin{align*}
\bE[M_{k+1}^{4}| \overline{\cF}_{k}, N_{T}= \ell]
& \le
C
\frac{\Delta b_{k}^{4}}{\Delta T_{k+1}^{2}} \frac{1}{\sigma(\theta_{k},T_{k+1})^{4}}
\le
C
\frac{\Delta T_{k}^{4} }{\Delta T_{k+1}^{2}} 
\frac{1 }{\sigma(\theta_{k},T_{k+1})^{4}}
\, , 
\\
\bE[V_{k+1}^{4}| \overline{\cF}_{k}, N_{T}= \ell]
& \le
C \frac{\sigma(\theta_{k-1},T_{k})^{8}}{\sigma(\theta_{k},T_{k+1})^{8}}
\frac{1}{\Delta T_{k+1}^{4}}
\, .
\end{align*}
By Assumption \ref{Assume:Finite Variance} and the fact that $\sigma$ is bounded from below implies that the $V$ term dominates the $M$ term, hence,
we obtain,
\begin{align}
\label{Eq:Full P Bound}
\bE[P_{k+1}^{4}| \overline{\cF}_{k}, N_{T}= \ell]
& \le 
C \frac{1}{f(\Delta T_{k})^{4}} 
\frac{\sigma(\theta_{k-1},T_{k})^{8}}{\sigma(\theta_{k},T_{k+1})^{8}}
\frac{1}{\Delta T_{k+1}^{4}}.
\end{align}
We are now able to consider bounding the term we originally set out to. Using the bound we obtained for $\beta^{2}$,
\begin{align}
& \bE \left[   \beta^{2} \prod_{k=2}^{N_{T}} P_{k}^{2} \Bigg | N_{T}=\ell \right] 
\notag
\\
& \le
\bE \Bigg[ \dfrac{
	C
}{
f(\Delta T_{N_{T}})^{2}} 
\Big( \Delta b_{N_{T}}^{2}
+
\frac{\sigma(\theta_{N_{T}-1},T_{N_{T}})^{4}}{\sigma(\theta_{N_{T}},T_{N_{T}+1})^{2}}
\frac{\min \left[
	1
	,
	\sigma(\theta_{N_{T}},T_{N_{T}+1})^{2} \Delta T_{N_{T}+1}
	\right]}{\Delta T_{N_{T}+1}}
\Big)
\prod_{k=2}^{N_{T}}
P_{k}^{2}
\Bigg | N_{T}=\ell \Bigg]
\label{Eq:Estimator with beta bound} 
\, .
\end{align}
One can view this product as having two components, one which does not depend on $\Delta T_{N_{T}+1}$ which comes from the $\Delta b_{N_{T}}$ and a component that does depend on $\Delta T_{N_{T}+1}$. In order to show that the second moment is finite we split these two components and show each of them is finite.






\emph{Step 3: Bounding each product in \eqref{Eq:Estimator with beta bound}}. 
Let us start by considering the product from the $\Delta b_{N_{T}}$ term
\begin{align*}
\bE \Bigg[ \dfrac{
	\Delta b_{N_{T}}^{2}
}{
f(\Delta T_{N_{T}})^{2}} 
\prod_{k=2}^{N_{T}}
P_{k}^{2}
\Bigg | N_{T}=\ell \Bigg] 
& =
\bE \Bigg[ \dfrac{
	1
}{
f(\Delta T_{N_{T}})^{2}} 
\bE[
\Delta b_{N_{T}}^{2} P_{N_{T}}^{2}
| \overline{\cF}_{N_{T}-1}, N_{T}= \ell]
\prod_{k=2}^{N_{T}-1}
P_{k}^{2}
\Bigg | N_{T}=\ell \Bigg] .
\end{align*}
Applying Cauchy-Schwarz to the internal expectation and using the previous bounds we obtain,
\begin{align}
\label{Eq:Bound b P term}
\bE[
\Delta b_{N_{T}}^{2} P_{N_{T}}^{2}
| \overline{\cF}_{N_{T}-1}, N_{T}= \ell]
& \le
C \frac{1}{f(\Delta T_{N_{T}-1})^{2}} 
\frac{\sigma(\theta_{N_{T}-2},T_{N_{T}-1})^{4}}{\sigma(\theta_{N_{T}-1},T_{N_{T}})^{4}}
\, .
\end{align}
Note that this bound and \eqref{Eq:Full P Bound} have no dependence on the Brownian motion, therefore we can isolate each $P_{k}$ by recursively conditioning, i.e.
\begin{align*}
&\bE \Bigg[ \dfrac{
	C
}{
f(\Delta T_{N_{T}})^{2}} 
\bE[
\Delta b_{N_{T}}^{2} P_{N_{T}}^{2}
| \overline{\cF}_{N_{T}-1}, N_{T}= \ell]
\prod_{k=2}^{N_{T}-1}
P_{k}^{2}
\Bigg | N_{T}=\ell \Bigg] 
\\
&
\le
\bE \Bigg[ \dfrac{
	C
}{
f(\Delta T_{N_{T}})^{2}} 
\frac{1}{f(\Delta T_{N_{T}-1})^{2}} 
\frac{\sigma(\theta_{N_{T}-2},T_{N_{T}-1})^{4}}{\sigma(\theta_{N_{T}-1},T_{N_{T}})^{4}}
\prod_{k=2}^{N_{T}-1}
P_{k}^{2}
\Bigg | N_{T}=\ell \Bigg] 
\\
&
=
\bE \Bigg[ \dfrac{
	C
}{
f(\Delta T_{N_{T}})^{2}} 
\frac{1}{f(\Delta T_{N_{T}-1})^{2}} 
\frac{\sigma(\theta_{N_{T}-2},T_{N_{T}-1})^{4}}{\sigma(\theta_{N_{T}-1},T_{N_{T}})^{4}}
\bE[
P_{N_{T}-1}^{2}
| \overline{\cF}_{N_{T}-2}, \Delta T_{N_{T}}, N_{T}= \ell]
\prod_{k=2}^{N_{T}-2}
P_{k}^{2}
\Bigg | N_{T}=\ell \Bigg] .
\end{align*}
Using our results and noting that most of the $\sigma$ terms cancel yields the following bound,
\begin{align*}
& \bE \Bigg[ \dfrac{
	\Delta b_{N_{T}}^{2}
}{
f(\Delta T_{N_{T}})^{2}} 
\prod_{k=2}^{N_{T}}
P_{k}^{2}
\Bigg | N_{T}=\ell \Bigg] 
\\
& \le
\bE \Bigg[
\dfrac{
	C^{N_{T}}
}{
f(\Delta T_{N_{T}})^{2}} 
\frac{1}{f(\Delta T_{N_{T}-1})^{2}} 
\frac{\sigma(\theta_{N_{T}-2},T_{N_{T}-1})^{4}}{\sigma(\theta_{N_{T}-1},T_{N_{T}})^{4}}
\prod_{k=2}^{N_{T}-1}
\frac{1}{f(\Delta T_{k-1})^{2}} 
\frac{\sigma(\theta_{k-2},T_{k-1})^{4}}{\sigma(\theta_{k-1},T_{k})^{4}}
\frac{1}{\Delta T_{k}^{2}}
\Bigg | N_{T}=\ell \Bigg] 
\\
&
=
\bE \Bigg[
\dfrac{
	C^{N_{T}}
}{
f(\Delta T_{N_{T}})^{2}} 
\frac{1}{f(\Delta T_{1})^{2}} 
\frac{\sigma(\theta_{0},T_{1})^{4}}{\sigma(\theta_{N_{T}-1},T_{N_{T}})^{4}}
\prod_{k=2}^{N_{T}-1}
\frac{1}{f(\Delta T_{k})^{2}} 
\frac{1}{\Delta T_{k}^{2}}
\Bigg | N_{T}=\ell \Bigg] \, .
\end{align*}
Recall the goal here is to ultimately bound this by a term of the form $C^{N_{T}}$, which holds provided all $\Delta T_{k}$ dependence is to a positive power. Recall that since $f$ is the density for the Gamma distribution with shape $\kappa$, we have that,
\begin{align*}
f(\Delta T_{k}) \ge C \Delta T_{k}^{\kappa -1}
\quad
\implies
\quad
\frac{1}{f(\Delta T_{N_{T}})^{2}} 
\le
C \Delta T_{N_{T}}^{2 - 2 \kappa} \, .
\end{align*}
Using the representation for $\sigma$ we obtain terms of the form $\Delta T_{k}^{2-2 \kappa -2 -4n}$, hence we require $ 2 \kappa -4n \ge 0$, which suggests $n \le -\kappa/2$. Since Assumption \ref{Assume:Finite Variance} implies these conditions on $n$ and $\kappa$ hold\footnote{Note that $\kappa = 1/2$ also implies $1/ f(\Delta T_{1}) \le C$.}, one obtains
\begin{align}
\label{Eq:Bound for b term}
& \bE \Bigg[ \dfrac{
	\Delta b_{N_{T}}^{2}
}{
f(\Delta T_{N_{T}})^{2}} 
\prod_{k=2}^{N_{T}}
P_{k}^{2}
\Bigg | N_{T}=\ell \Bigg] 
\le
\bE \Big[
C^{N_{T}}
\Big | N_{T}=\ell \Big] \, .
\end{align}
Showing this is finite is done in \cite{DoumbiaEtAl2017}. As it turns out the other term in \eqref{Eq:Estimator with beta bound} also dominates this term, hence we do not discuss it further.

For the second term in \eqref{Eq:Estimator with beta bound} we note that the $\sigma$ terms do not depend on the Brownian motion, hence we can again condition to isolate the various $P_{k}$ terms, hence, 
\begin{align}
& \bE \Bigg[ \dfrac{
	C
}{
f(\Delta T_{N_{T}})^{2}}
\frac{\sigma(\theta_{N_{T}-1},T_{N_{T}})^{4}}{\sigma(\theta_{N_{T}},T_{N_{T}+1})^{2}}
\frac{\min \left[
	1
	,
	\sigma(\theta_{N_{T}},T_{N_{T}+1})^{2} \Delta T_{N_{T}+1}
	\right]}{\Delta T_{N_{T}+1}}
\prod_{k=2}^{N_{T}}
P_{k}^{2}
\Bigg | N_{T}=\ell \Bigg] 
\notag
\\
&
\le
\bE \Bigg[ \frac{
	C
}{
f(\Delta T_{N_{T}})^{2}}
\frac{\sigma(\theta_{N_{T}-1},T_{N_{T}})^{4}}{\sigma(\theta_{N_{T}},T_{N_{T}+1})^{2}}
\frac{\min \left[
	1
	,
	\sigma(\theta_{N_{T}},T_{N_{T}+1})^{2} \Delta T_{N_{T}+1}
	\right]}{\Delta T_{N_{T}+1}}
\notag
\\
&
\qquad \qquad
\times
\prod_{k=2}^{N_{T}}
\frac{1}{f(\Delta T_{k-1})^{2}} 
\frac{\sigma(\theta_{k-2},T_{k-1})^{4}}{\sigma(\theta_{k-1},T_{k})^{4}}
\frac{1}{\Delta T_{k}^{2}}
\Bigg | N_{T}=\ell \Bigg] 
.
\label{Eq:sigma term}
\end{align}
By cancelling repeating $\sigma$ terms in the product and again using $1/ f(\Delta T_{1}) \le C$, we obtain the following simpler result,
\begin{align}
\label{Eq:Bounding Over the Product}
\eqref{Eq:sigma term}
\le
C \bE \left[ 
\frac{\sigma(\theta_{0},T_{1})^{4}}{\sigma(\theta_{N_{T}},T_{N_{T}+1})^{2}}\frac{\min \left[
	1
	,
	\sigma(\theta_{N_{T}},T_{N_{T}+1})^{2} \Delta T_{N_{T}+1}
	\right]}{\Delta T_{N_{T}+1}}
\prod_{k=2}^{N_{T}}
\frac{1}{f(\Delta T_{k})^{2}} 
\frac{1}{\Delta T_{k}^{2}} \Bigg | N_{T}=\ell \right] 
\, . 
\end{align}
Using the fact that $\sigma(\theta_{0},T_{1})=\sigma_{0}$ and $f$ is the density for the Gamma distribution we can bound \eqref{Eq:Bounding Over the Product} by, 
\begin{align}
\label{Eq:Full bound with nu}
&
\bE \left[ C^{N_{T}}
\frac{\min \left[
	1
	,
	\sigma(\theta_{N_{T}},T_{N_{T}+1})^{2} \Delta T_{N_{T}+1}
	\right]}
{\sigma(\theta_{N_{T}},T_{N_{T}+1})^{2}} \Delta T_{N_{T}+1}^{-1}
\prod_{k=2}^{N_{T}}
\Delta T_{k}^{-2 \kappa}
\Bigg | N_{T}=\ell \right] 
\notag
\\
&
\le
\bE \left[ C^{N_{T}}
\frac{\sigma(\theta_{N_{T}},T_{N_{T}+1})^{\nu} \Delta T_{N_{T}+1}^{\nu /2}}
{\sigma(\theta_{N_{T}},T_{N_{T}+1})^{2}} \Delta T_{N_{T}+1}^{-1}
\prod_{k=2}^{N_{T}}
\Delta T_{k}^{-2 \kappa}
\Bigg | N_{T}=\ell \right] 
\quad \text{for } \nu \in [0,2] \, ,
\end{align}
where the inequality comes from the observation that,
\begin{align*}
\min \left[
1
,
\sigma(\theta_{N_{T}},T_{N_{T}+1})^{2} \Delta T_{N_{T}+1}
\right]
\le 
\sigma(\theta_{N_{T}},T_{N_{T}+1})^{\nu} \Delta T_{N_{T}+1}^{\nu /2} \quad \text{for any } \nu \in [0,2] \, .
\end{align*}
The presence of $\Delta T_{N_{T}+1}^{-1}$ makes \eqref{Eq:Full bound with nu} more challenging. Of course, one could take $\nu=2$ to remove $\Delta T_{N_{T}+1}^{-1}$, however, this also removes $\sigma$ and since $\kappa >0$ we are still left with an unbounded product. Therefore we must chose $\nu$ carefully and apply a delicate argument to appropriately bound \eqref{Eq:Full bound with nu}. 

One can note the similarity between \eqref{Eq:Full bound with nu} and \eqref{Eq:Bound for b term}. However, \eqref{Eq:Full bound with nu} is more complex and as it turns out, the bound we eventually achieve for it dominates \eqref{Eq:Bound for b term}. We therefore complete the proof showing \eqref{Eq:Full bound with nu} is bounded, since this implies \eqref{Eq:Bound for b term} is bounded.






\emph{Step 4: Interval splitting.}
Recall we are interested in proving convergence of the sum
\begin{align*}
\sum_{\ell =1}^{\infty} \bE \left[ \left( \beta \prod_{k=2}^{\ell} P_{k} \right)^{2} \Bigg | N_{T}=\ell \right] \bP[N_{T}=\ell] \, .
\end{align*}
Let us split this into two components, $\ell =1$ and $\ell \ge 2$. When $\ell =1$ we obtain nothing from the product and are thus only showing that $\beta$ is square integrable, such is obvious from our previous calculations. We now concentrate on the case $\ell \ge 2$. Recall that for $i=1, \dots, M$, if $Y_{i} \sim \Gamma (a,b)$ i.i.d. then $\sum_{i=1}^{M} Y_{i} \sim \Gamma (aM,b)$ and fix $\ell \ge 2$, we can then partition the expectation as follows,
\begin{align*}
& \bE \left[ \left( \beta \prod_{k=2}^{\ell} P_{k} \right)^{2} \Bigg | N_{T}=\ell \right]
=
\bE \left[ \left( \beta \prod_{k=2}^{\ell} P_{k} \right)^{2} \Bigg | N_{T}=\ell , \Delta T_{N_{T}+1} \ge \frac{T}{\ell}\right]
\bP\left[ \Delta T_{N_{T}+1} \ge \frac{T}{\ell} \Big | N_{T}=\ell \right]
\\ & \quad
+
\sum_{m=1}^{\infty}
\bE \left[ \left( \beta \prod_{k=2}^{\ell} P_{k} \right)^{2} \Bigg | N_{T}=\ell , \frac{T}{\ell^{m+1}} \le \Delta T_{N_{T}+1} < \frac{T}{\ell^{m}}\right]
\bP\left[ \frac{T}{\ell^{m+1}} \le \Delta T_{N_{T}+1} < \frac{T}{\ell^{m}} \Big | N_{T}=\ell \right]
\, .
\end{align*}
Firstly, we note that when $\Delta T_{N_{T}+1} \ge T/ \ell$, the expectation is simple to bound since we can take the minimum as $1$ (the $\nu=0$ case in \eqref{Eq:Full bound with nu}) then use the fact $\sigma(\theta_{N_{T}},T_{N_{T}+1})^{-2}= \sigma_{0}^{-2}\prod_{i=1}^{\ell} \Delta T_{i}^{2n}$ and $\kappa < -n$ by Assumption \ref{Assume:Finite Variance}. Hence the following bound holds,
\begin{align*}
\bE \left[ \left( \beta \prod_{k=2}^{\ell} P_{k} \right)^{2} \Bigg | N_{T}=\ell , \Delta T_{N_{T}+1} \ge \frac{T}{\ell}\right]
\bP\left[ \Delta T_{N_{T}+1} \ge \frac{T}{\ell} \Big | N_{T}=\ell \right]
\le
\ell C^{\ell}
\, .
\end{align*}
For the case $m \ge 1$, we have that
\begin{align*}
\bP\left[ \frac{T}{\ell^{m+1}} \le \Delta T_{N_{T}+1} < \frac{T}{\ell^{m}} \Big | N_{T}=\ell  \right]
=
\bP\left[  T-\frac{T}{\ell^{m}} \le \sum_{i=1}^{\ell} \Delta T_{i} <T- \frac{T}{\ell^{m+1}} \Big | N_{T}=\ell  \right] \, .
\end{align*}
Due to the fact $\kappa = 1/2$ by Assumption \ref{Assume:Finite Variance}, the distribution of $\sum_{i=1}^{\ell} \Delta T_{i} $ is Gamma with shape parameter at least $1$, therefore the density has a finite maximum, unfortunately the conditioning makes this probability difficult to deal with. We therefore expand,
\begin{align*}
& \bP\left[  T-\frac{T}{\ell^{m}} \le \sum_{i=1}^{\ell} \Delta T_{i} <T- \frac{T}{\ell^{m+1}} \Big | N_{T}=\ell  \right]
\\ &
\qquad =
\frac{1}{\bP[N_{T}=\ell]}
\bP\left[  T-\frac{T}{\ell^{m}} \le \sum_{i=1}^{\ell} \Delta T_{i} <T- \frac{T}{\ell^{m+1}} , ~
\sum_{i=1}^{\ell} \Delta T_{i} <T, ~
\sum_{i=1}^{\ell+1} \Delta T_{i} \ge T
\right] 
\\ &
\qquad \le
\frac{1}{\bP[N_{T}=\ell]}
\bP\left[  T-\frac{T}{\ell^{m}} \le \sum_{i=1}^{\ell} \Delta T_{i} <T- \frac{T}{\ell^{m+1}} 
\right]\, .
\end{align*}
Using this form we have removed the conditional dependence on the number of jumps and therefore we can use the distribution of $\sum_{i=1}^{\ell} \Delta T_{i} $. We note that for $\ell$ large the density of the distribution at point $T$ will be larger than values less than $T$, further, since the density has a finite maximum, for $\ell$ smaller we can bound by some constant multiplied by the value at point $T$, thus,
\begin{align*}
\bP\left[  T-\frac{T}{\ell^{m}} \le \sum_{i=1}^{\ell} \Delta T_{i} <T- \frac{T}{\ell^{m+1}} 
\right]
\le
C \ell^{-m} f(T)
\le
C \ell^{-m} \frac{T^{\ell \kappa -1} e^{-T/\eta}}{\eta^{\ell \kappa} \Gamma(\ell \kappa)}\, ,
\end{align*}
where we have used the p.d.f. of a Gamma random variable to obtain the last inequality.
Similar to the case $\ell =1$ we can bound the expectation by
\begin{align*}
& \bE \left[ \left( \beta \prod_{k=2}^{N_{T}} P_{k} \right)^{2} \Bigg | N_{T}=\ell , \frac{T}{\ell^{m+1}} \le \Delta T_{N_{T}+1} < \frac{T}{\ell^{m}}\right]
\\ &
\qquad
\le
\bE \left[ C^{N_{T}}
\Delta T_{N_{T}+1}^{-1+\nu/2}
\prod_{k=2}^{N_{T}}
\Delta T_{k}^{-(2- \nu)n-2 \kappa}
\Bigg | N_{T}=\ell  ,\frac{T}{\ell^{m+1}} \le \Delta T_{N_{T}+1} < \frac{T}{\ell^{m}} \right] 
\, .
\end{align*}
A simple requirement for the product to be bounded is $-(2- \nu)n-2 \kappa \ge 0$, by Assumption \ref{Assume:Finite Variance} $\kappa =1/2$, hence $-n \ge 1/(2 - \nu)$. As it turns out, taking $\nu =1$ is useful to complete the proof, therefore we require $n \le -1$, which holds by Assumption \ref{Assume:Finite Variance}. This set of $\kappa$, $\nu$ and $n$ also allow us to bound \eqref{Eq:Bound for b term}, hence we only considered \eqref{Eq:Full bound with nu}.

The only term we have to consider in the expectation is $\Delta T_{N_{T}+1}^{-1+\nu/2}$, but by our conditioning this is bounded by $T \ell^{(1-\nu/2)(m+1)}$, hence for fixed $\ell \ge 2$ and letting $\nu=1$ we obtain the following,
\begin{align*}
& \bE \left[ \left( \beta \prod_{k=2}^{N_{T}} P_{k} \right)^{2} \Bigg | N_{T}=\ell \right]
\le
C^{\ell} \ell
+
\frac{1}{\bP[N_{T}=\ell]}
\sum_{m=1}^{\infty}
C^{\ell} \ell^{(1/2)(m+1)} \ell^{-m}
\frac{T^{\ell \kappa -1} e^{-T/\eta}}{\eta^{\ell \kappa} \Gamma(\ell \kappa)}
\, .
\end{align*}
One can easily see that the sum in $m$ converges since $(1/2)(m+1)-m \le 0$ for $m \ge 1$ and $\ell \ge 2$, the sum can be easily bounded by $\sum_{m=1}^{\infty}2^{-(1/2)m +1/2}=C$ for any $\ell \ge 2$. 
One can compare this to the result in \cite{DoumbiaEtAl2017}*{Proposition 4.1} where the authors obtain a bound of the form $C^{\ell}$, hence our bound is not as strong but it is still good enough to ensure convergence.






\emph{Step 5: The sum over $N_{T}$ converges}.
The final step of the proof is to show that the overall sum converges. We proceed by observing the following (see \cite{DoumbiaEtAl2017}*{Proposition 4.1}),
\begin{align*}
\bP[N_{T}=\ell] \le \frac{C^{\ell \kappa}}{\ell \kappa \Gamma(\ell \kappa)} \, .
\end{align*}
Using a generalisation of Stirling's formula one can approximate $\Gamma(z) \sim z^{z-1/2} e^{-z}\sqrt{2 \pi}$. Hence we can bound
\begin{align*}
\bE \left[\beta \prod_{k=2}^{N_{T}} P_{k} \1_{\{N_{T} \ge 1\}} \right] 
&
\le
\sum_{\ell =1}^{\infty} C^{\ell} \ell
\frac{C^{\ell \kappa}}{\ell \kappa \Gamma(\ell \kappa)}
+
\frac{\bP[N_{T}=\ell]}{\bP[N_{T}=\ell]}
\sum_{m=1}^{\infty}
C^{\ell} \ell^{(1/2)(m+1)} \ell^{-m}
\frac{T^{\ell \kappa -1} e^{-T/\eta}}{\eta^{\ell \kappa} \Gamma(\ell \kappa)}
\\
&
\le
\sum_{\ell =1}^{\infty} C^{\ell} 
\frac{C^{\ell \kappa}}{ \kappa \Gamma(\ell \kappa)}
\, ,
\end{align*}
and using Stirling's formula,
\begin{align*}
C^{\ell}
\frac{C^{\ell \kappa}}{\kappa \Gamma(\ell \kappa)} 
\sim
C^{\ell}
\frac{C^{\ell \kappa} e^{\ell \kappa}}{\kappa (\ell \kappa)^{\ell \kappa -1/2} \sqrt{2 \pi}} 
\le
\left(
\frac{C^{1/\kappa}  e^{1}}{\ell \kappa}
\right)^{\ell \kappa -1/2}
C^{1/(2\kappa)}e^{1/2}
\, ,
\end{align*} 
since $\kappa=1/2$ this gives a sequence that converges under summation.
\end{proof}

\begin{remark}
	[Optimal $\sigma_{0}$]
	One can see from the variance calculations that the $\frac{\sigma(\theta,T_{0})^{4}}{\sigma(\theta,T_{N_{T}})^{2}}$ will leave a $\sigma_{0}^{2}$ term behind. Thus as one would expect the variance will be minimised by taking $\sigma_{0}$ smaller, however, to deal with terms involving nonlinearities in $\partial_{x} v$ one obtains terms of the form $\frac{1}{\sigma}$ thus an optimisation needs to be performed in order to set $\sigma_{0}$ at the correct level. Crucially however, the expected value (bias) is not effected by this choice.
\end{remark}


\subsection{Estimator solves the PDE under enough integrability}

At this point we have only proved that the estimator can be approximated via Monte Carlo. We now show that given some extra integrability conditions the estimator solves PDE \eqref{Eq:Transport PDE}. The final step is to show the said integrability conditions hold.


Theorem \ref{theo:StochasticRepresentation} is the analogous result to Theorem 3.5 in \cite{LabordereEtAl2016}, however, the representation we derive below is more complex. The reason for the added complexity is the antithetic as well as the control variate on the final jump. Where as the control variate keeps the final Malliavin weight the same, the antithetic changes the weight, this then requires us to have extra terms that \cite{LabordereEtAl2016} does not have.
\begin{theorem}
	\label{theo:StochasticRepresentation}
	Let Assumptions \ref{Assume:Main SDE Assumptions}, \ref{Assume:Bounded g} and \ref{Assume:Finite Variance} hold. Define the following random variables,
	\begin{align*}
	\tilde{\psi}^{t,x}
	:=
	&
	\Bigg(
	\frac{\Delta g_{T_{N_{T}+1}}}{2\overline{F}(\Delta T_{N_{T}+1})}
	\frac{\Delta b_{N_{T}} \cW_{N_{T}+1}^{1} - \frac{1}{2} \sigma(\theta_{N_{T}-1},T_{N_{T}})^{2} \cW_{N_{T}+1}^{2}}{f(\Delta T_{N_{T}})} 
	\\
	&
	+ \frac{\Delta \hat{g}_{T_{N_{T}+1}}}{2\overline{F}(\Delta T_{N_{T}+1})}
	\frac{-\Delta b_{N_{T}} \cW_{N_{T}+1}^{1} - \frac{1}{2} \sigma(\theta_{N_{T}-1},T_{N_{T}})^{2} \cW_{N_{T}+1}^{2}}{f(\Delta T_{N_{T}})}
	\Bigg)
	\prod_{k=2}^{N_{T}}
	\frac{\Delta b_{k-1} \cW_{k}^{1} - \frac{1}{2} \sigma(\theta_{k-2},T_{k-1})^{2} \cW_{k}^{2}}{f(\Delta T_{k-1})}
	\, ,
	\\
	\psi^{t,x}
	:=
	&
	\1_{\{N_{T} =0 \}}
	\frac{g(\bar{X}_{T_{N_{T}+1}})}{\overline{F}(\Delta T_{N_{T}+1})}
	+
	\1_{\{N_{T} \ge 1\}}
	\tilde{\psi}^{t,x}
	\, ,
	\end{align*}
	and
	\begin{align*}
	\Phi_{1}^{T_{N_{T}}, \bar{X}_{T_{N_{T}}}} = \frac{\Delta g_{T_{N_{T}+1}}-\Delta \hat{g}_{T_{N_{T}+1}}}{2\overline{F}(\Delta T_{N_{T}+1})}
	\quad
	\text{and}
	\quad
	\Phi_{2}^{T_{N_{T}}, \bar{X}_{T_{N_{T}}}} = \frac{\Delta g_{T_{N_{T}+1}}+\Delta \hat{g}_{T_{N_{T}+1}}}{2\overline{F}(\Delta T_{N_{T}+1})}
	\end{align*}
	where
	\begin{align*}
	\Delta g_{T_{N_{T}+1}}
	:=
	g(\bar{X}_{T_{N_{T}+1}})-g(\bar{X}_{T_{N_{T}}} + b(T_{N_{T}}) \Delta T_{N_{T}+1}) \, ,
	\\
	\Delta \hat{g}_{T_{N_{T}+1}}
	:=
	g(\hat{X}_{T_{N_{T}+1}})-g(\bar{X}_{T_{N_{T}}} + b(T_{N_{T}}) \Delta T_{N_{T}+1}) \, ,
	\end{align*}
	the first and second order Malliavin weights are given by,
	\begin{align}
	\label{Eq:First and Second Malliavin Weights}
	\cW_{k+1}^{1}
	=
	\sigma(\theta_{k},T_{k+1})^{-1} \frac{\Delta W_{T_{k+1}}}{ \Delta T_{k+1}} 
	\quad
	\text{and}
	\quad
	\cW_{k+1}^{2}
	=
	\sigma(\theta_{k},T_{k+1})^{-2}
	\left(
	\frac{\Delta W_{T_{k+1}}^{2} - \Delta T_{k+1}}{\Delta T_{k+1}^{2}}
	\right) \, .
	\end{align}
	The superscript in $\psi$, $\tilde{\psi}$, $\Phi_{1}$ and $\Phi_{2}$ denotes the initial condition for the SDE, $\bar{X}$. Further assume that,
	\begin{align*}
	& \psi^{t,x}, 
	~ ~\tilde{\psi}^{t,x} \cW_{1}^{1} , 
	~ ~\tilde{\psi}^{t,x} \cW_{1}^{2} ,
	~ ~f(\Delta T_{1})^{-1}\Delta b_{1} \tilde{\psi}^{T_{1},\bar{X}_{T_{1}}} \cW_{2}^{1} ,
	~ ~f(\Delta T_{1})^{-1}\sigma(\theta_{0}, T_{1})^{2}  \tilde{\psi}^{T_{1},\bar{X}_{T_{1}}} \cW_{2}^{2} ,
	\\
	&
	 \Phi_{1}^{T_{N_{T}},\bar{X}_{T_{N_{T}}}} \cW_{N_{T}+1}^{1} ,
	~ ~ \Phi_{2}^{T_{N_{T}},\bar{X}_{T_{N_{T}}}} \cW_{N_{T}+1}^{2} \, ,
	\end{align*}
	are uniformly integrable and that $\psi^{T_{1}, \bar{X}_{T_{1}}}$, $\Delta b_{2} \tilde{\psi}^{T_{2},\bar{X}_{T_{2}}} \cW_{3}^{1}$,
	$\sigma(\theta_{1}, T_{2})^{2}  \tilde{\psi}^{T_{2},\bar{X}_{T_{2}}} \cW_{3}^{2}$ are $\bP$-a.s. uniformly integrable and $\tilde{\psi}^{T_{1},\bar{X}_{T_{1}}} \cW_{2}^{1}$ and $\tilde{\psi}^{T_{1},\bar{X}_{T_{1}}} \cW_{2}^{2}$ are $\bP$-a.s. integrable.
	
	Then, the function $\hat{v}(t,x):= \bE[\psi^{t,x} | \cF_{t}]$ solves the PDE \eqref{Eq:High Dimension Transport PDE}.
\end{theorem}

\begin{remark}[$\bP$-a.s. (uniformly) integrable]
	Note that some of the processes stated in the theorem, for example $\psi^{T_{1}, \bar{X}_{T_{1}}}$ and $\tilde{\psi}^{T_{1},\bar{X}_{T_{1}}} \cW_{2}^{2}$ depend on random ``initial conditions''. Hence some of these processes are unbounded, but are finite up to a null set. For example, when we state $\tilde{\psi}^{T_{1},\bar{X}_{T_{1}}} \cW_{2}^{2}$ is $\bP$-a.s. integrable, we mean that, $\bE[|\tilde{\psi}^{T_{1},\bar{X}_{T_{1}}} \cW_{2}^{2}| \,  | \cF_{T_{1}}] < \infty$ $\bP$-a.s. and similar for the uniform integrability condition. Recall that $\bP$ is the product measure $\bP_{W} \otimes \bP_{f}$.	
\end{remark}

This theorem only shows that the estimator gives rise to the solution of the PDE under certain integrability assumptions. In order to finish our proof we need to show that such integrability conditions hold (Theorem \ref{Thm:Integrability Hold}).
Although it is $\psi$ that solves the PDE, our proof relies on various intermediary steps requiring additional integrability on $\psi \cW$. Since one does not have this in general, we introduce the seemingly arbitrary $\tilde{\psi}$ and $\Phi$ which have the required integrability. Therefore, throughout the proof we show that one can view these additional processes as $\psi \cW$ with a control variate and perform the various steps on $\tilde{\psi}$ and $\Phi$.

\begin{remark}
	The Malliavin weights are given by \eqref{Eq:Explicit Malliavin Weights} since our unbiased estimation puts us in the simple setting where the SDE has constant coefficients (see \cite{FournieEtAl1999}).
\end{remark}

\begin{proof}
	The main idea of this proof is to first show a stochastic representation for the PDE, then show that this representation and $\bE[\psi^{t,x}| \cF_{t}]$ are equivalent. Following Section \ref{Sec:Heuristics}, since a $C_{b}^{1,3}$ solution is assumed to exist, one can take constants $b_{0}$ and $\sigma_{0}$ and define the following PDE (equivalent to \eqref{Eq:High Dimension Transport PDE}),
	\begin{equation*}
	\begin{cases}
	\partial_{t} v(t,x) + b_{0} \partial_{x}v(t,x) + \frac{1}{2}\sigma_{0}^{2} \partial_{xx}v(t,x) +
	(b(t)-b_{0})\partial_{x}v(t,x) - \frac{1}{2} \sigma_{0}^{2} \partial_{xx}v(t,x)
	=0  
	\, ,
	\\
	v(T,x)=g(x) \, .
	\end{cases}
	\end{equation*}
	Assume that these constants $b_{0}$ and $\sigma_{0}$ are adapted to the filtration $\cF_{t}$ (as defined at the start of Section \ref{Sec:Main Results}). Define $\tilde{X}$ as the solution to the SDE on $s \in [t,T]$
	\begin{align*}
	\dd \tilde{X}_{s} = b_{0} \dd s + \sigma_{0} \dd W_{s} \, , \quad \tilde{X}_{t}=x \, .
	\end{align*}	
	again since $v \in C_{b}^{1,3}$, one obtains from the Feynman-Kac formula,
	\begin{align*}
	v(t,x)=
	\bE_{W} \left[
	g(\tilde{X}_{T}) 
	+
	\int_{t}^{T} 
	(b(s)-b_{0})\partial_{x}v(s,\tilde{X}_{s}) - \frac{1}{2} \sigma_{0}^{2} \partial_{xx}v(s,\tilde{X}_{s}) \dd s 
	~
	\Big | \cF_{t}
	\right] \, .
	\end{align*}
	It is important to note that we have not assigned values to the constants $b_{0}$ and $\sigma_{0}$ here, only that they are adapted to the initial filtration.
	Using standard branching arguments, we introduce a random variable independent of Brownian motion, corresponding to the life of the particle which allows us to rewrite the previous expression as\footnote{Where $\bE$ is the expectation in the product space of the two random variables.},
	\begin{align}
	\label{Eq:Forward Representation}
	v(t,x)=
	\bE \left[
	\frac{g(\tilde{X}_{T})}{\overline{F}(\Delta T_{1})} \1_{\{T_{1}=T\}} 
	+
	\frac{\1_{\{T_{1}<T\}}}{f(\Delta T_{1})}
	\left \{
	(b(T_{1})-b_{0})\partial_{x}v(T_{1},\tilde{X}_{T_{1}}) - \frac{1}{2} \sigma_{0}^{2} \partial_{xx}v(T_{1},\tilde{X}_{T_{1}}) 
	\right \} 
	~
	\Big | \cF_{t}
	\right] \, .
	\end{align}
	As before, the representation does not depend on the value of the constants, therefore let us take $b_{0}:= b(t)$ and $\sigma_{0}:= \sigma_{0}$ (in the sense of \eqref{Eq:Sigma Defn}), thus $\tilde{X}$ is equivalent to $\bar{X}$.
	
	This can be thought of as the forward representation, the goal now is to reach the same representation going backwards. Namely, starting from the estimator $\psi^{t,x}$, we want to remove the Malliavin weights and obtain the same relationship. We break the remainder of the proof into several steps.

	\emph{Step 1: Continuity of the functions.}
	We start by noting that between any two mesh points, the SDE is continuous w.r.t. its initial condition $(T_{k}, \bar{X}_{T_{k}})$, which is clear from the fact that it is just an SDE with constant coefficients. This along with the uniform integrability assumption of $\psi$ implies that the function $\hat{v}$ is jointly continuous. This stems from the fact that we can define $\psi_{n}^{t,x}$ as $\psi^{t,x}$ but with the $N_{T}$ replaced by $N_{T} \wedge n$, hence $\psi^{t,x}= \lim_{n \rightarrow \infty} \psi_{n}^{t,x}$. Then for each $n$ we have a finite product of jointly continuous functions, which is therefore jointly continuous. Then uniform integrability allows us to take the limit as $n \rightarrow \infty$ inside to conclude that $(t,x) \rightarrow \bE[\psi^{t,x} | \cF_{t}]$ must also be a jointly continuous function.
	
The weights $\cW^{i}$ for $i=1,2$ are also continuous w.r.t. the initial condition. Thus by arguing in a similar way to above we have $\bE[\tilde{\psi}^{t,x}\cW_{1}^{i} | \cF_{t}]$ and $\bE[\Phi_{i}^{T_{N_{T}}, \bar{X}_{T_{N_{T}}}} \cW_{N_{T}+1}^{i} | \cF_{t}]$ are jointly continuous by the uniform integrability assumption.

	\emph{Step 2: Rewriting the representation.}
	By construction of $\psi$, there are two main cases, either the particle goes through a regime switch, which implies $\{ N_{T} \ge 1 \}$ or it ``survives'' until the end, $\{N_{T}=0\}$. The key difference to the representation is the introduction of the variance reduction techniques when $\{N_{T} \ge 1\}$, this is also the distinction between $\psi$ and $\tilde{\psi}$. Hence the representation is,
	\begin{align*}
	\hat{v}(t,x)= 
	& \bE \Bigg[
	\1_{\{N_{T} =0 \}}
	\frac{g(\bar{X}_{T_{N_{T}+1}})}{\overline{F}(\Delta T_{N_{T}+1})}
	\\
	&
	\qquad +
	\1_{\{N_{T} \ge 1\}}
	\Bigg(
	\frac{\Delta g_{T_{N_{T}+1}}}{2\overline{F}(\Delta T_{N_{T}+1})}
	\frac{\Delta b_{N_{T}} \cW_{N_{T}+1}^{1} - \frac{1}{2} \sigma(\theta_{N_{T}-1},T_{N_{T}})^{2} \cW_{N_{T}+1}^{2}}{f(\Delta T_{N_{T}})}
	\\
	&
	\qquad \qquad \qquad \qquad 
	+ \frac{\Delta \hat{g}_{T_{N_{T}+1}}}{2\overline{F}(\Delta T_{N_{T}+1})}
	\frac{-\Delta b_{N_{T}} \cW_{N_{T}+1}^{1} - \frac{1}{2} \sigma(\theta_{N_{T}-1},T_{N_{T}})^{2} \cW_{N_{T}+1}^{2}}{f(\Delta T_{N_{T}})}
	\Bigg)
	\\
	&
	\qquad \qquad
	\times 
	\prod_{k=2}^{N_{T}}
	\frac{\Delta b_{k-1} \cW_{k}^{1} - \frac{1}{2} \sigma(\theta_{k-2},T_{k-1})^{2} \cW_{k}^{2}}{f(\Delta T_{k-1})}
	~
	\Bigg | ~
	X_{t}=x, ~ \sigma(\theta_{0},t)
	\Bigg] \, ,
	\end{align*}
	where we are using conditioning to state the initial condition of the SDE. In order to save space in the future we will stick to conditioning $\cF_{t}$.
	Concentrating on the case $\{N_{T} \ge 1\}$, then the random variable $\Delta T_{1}$ exists and satisfies $t < T_{1} < T$. Hence we can consider the filtration up to that point and by the tower property rewrite the $\{N_{T} \ge 1\}$ term in the expectation as,
	\begin{align}
	\bE \Bigg[
	\1_{\{ N_{T} \ge 1 \}} 
	\frac{1}{f(\Delta T_{1})}
	\Bigg \{ \hspace{1cm} & \notag
	\\
	\Delta b_{1} 
	\bE \Bigg[
	\1_{\{N_{T} =1 \}} &
	\frac{\Delta g_{T_{N_{T}+1}}-\Delta \hat{g}_{T_{N_{T}+1}}}{2\overline{F}(\Delta T_{N_{T}+1})} 
	\cW_{N_{T}+1}^{1}
	\notag
	\\
	+ \1_{\{N_{T} > 1\}}
	\Bigg( &
	\frac{\Delta g_{T_{N_{T}+1}}}{2\overline{F}(\Delta T_{N_{T}+1})}
	\frac{\Delta b_{N_{T}} \cW_{N_{T}+1}^{1} - \frac{1}{2} \sigma(\theta_{N_{T}-1},T_{N_{T}})^{2} \cW_{N_{T}+1}^{2}}{f(\Delta T_{N_{T}})}  \notag
	\\
	&
	+
	\frac{\Delta \hat{g}_{T_{N_{T}+1}}}{2\overline{F}(\Delta T_{N_{T}+1})}
	\frac{-\Delta b_{N_{T}} \cW_{N_{T}+1}^{1} - \frac{1}{2} \sigma(\theta_{N_{T}-1},T_{N_{T}})^{2} \cW_{N_{T}+1}^{2}}{f(\Delta T_{N_{T}})}
	\Bigg) \notag
	\\
	& ~
	\times \prod_{k=3}^{N_{T}}
	\frac{\Delta b_{k-1} \cW_{k}^{1} - \frac{1}{2} \sigma(\theta_{k-2},T_{k-1})^{2} \cW_{k}^{2}}{f(\Delta T_{k-1})}
	\cW_{2}^{1}
	~
	\Bigg | ~
	\cF_{T_{1}}
	\Bigg]  \notag
	\\
	-
	\frac{1}{2} \sigma(\theta_{0},T_{1})^{2} 
	\bE \Bigg[
	\1_{\{N_{T} =1 \}}&
	\frac{\Delta g_{T_{N_{T}+1}} + \Delta \hat{g}_{T_{N_{T}+1}}}{2\overline{F}(\Delta T_{N_{T}+1})} 
	\cW_{N_{T}+1}^{2}
	\notag
	\\
	+ \1_{\{N_{T} > 1\}}
	\Bigg( &
	\frac{\Delta g_{T_{N_{T}+1}}}{2\overline{F}(\Delta T_{N_{T}+1})} 
	\frac{\Delta b_{N_{T}} \cW_{N_{T}+1}^{1} - \frac{1}{2} \sigma(\theta_{N_{T}-1},T_{N_{T}})^{2} \cW_{N_{T}+1}^{2}}{f(\Delta T_{N_{T}})} 
	\notag
	\\
	&
	+
	\frac{\Delta \hat{g}_{T_{N_{T}+1}}}{2\overline{F}(\Delta T_{N_{T}+1})}
	\frac{-\Delta b_{N_{T}} \cW_{N_{T}+1}^{1} - \frac{1}{2} \sigma(\theta_{N_{T}-1},T_{N_{T}})^{2} \cW_{N_{T}+1}^{2}}{f(\Delta T_{N_{T}})}
	\Bigg) \notag
	\\
	& ~
	\times \prod_{k=3}^{N_{T}}
	\frac{\Delta b_{k-1} \cW_{k}^{1} - \frac{1}{2} \sigma(\theta_{k-2},T_{k-1})^{2} \cW_{k}^{2}}{f(\Delta T_{k-1})}
	\cW_{2}^{2}
	\Bigg | ~
	\cF_{T_{1}}
	\Bigg] 
	~ ~ ~
	\Bigg \}
	~
	\Bigg | ~
	\cF_{t}
	\Bigg]
	\label{Eq:At Least One Branching}
	\, ,
	\end{align} 
	where we have used that $\Delta b_{1}$ and $\sigma(\theta_{0},T_{1})$ are bounded and our integrability assumptions on $\Phi$ and $\tilde{\psi}^{T_{1}, \bar{X}_{T_{1}}}$ to apply the tower property.
	We see here that the antithetic variable is causing extra difficultly since we need to treat the case $N_{T}=1$ separately.

	\emph{Step 3: Existence and continuity of derivatives.}
	In order to obtain the required expression we must also understand the derivatives of the function, hence we must show these derivatives exist and obtain a representation for them.	 
	One can identify the terms inside the conditional expectations as $\Phi_{i}^{T_{N_{T}}, \bar{X}_{T_{N_{T}}}} \cW_{N_{T}+1}^{i}$ and $\tilde{\psi}^{T_{1}, \bar{X}_{T_{1}}} \cW_{2}^{i}$ for $i \in \{1,2\}$.
	
	Let us denote by $\eta(T_{1}, \bar{X}_{T_{1}}) := \bE[\psi^{t,x} | \cF_{T_{1}}]$, notice that for the same reasons $\psi^{t,x}$ is a continuous function of $x$, $\eta(T_{1}, \bar{X}_{T_{1}})$ is continuous w.r.t. $\bar{X}_{T_{1}}$ (which is in turn continuous w.r.t. $x$). Let us now consider derivatives of this function w.r.t. $x$.
	However, one should note that this expectation is on the product space of random variables $T_{i}$ and $W$. While the Malliavin automatic differentiation results only hold differentiating $\bE_{W}[ \cdot]$. Therefore we must swap the derivative with the expectation $\bE_{f}$, which we have proved to be valid (actually shown a more general case) in Lemma \ref{Lem:Swap derivative and integral} under the assumed integrability. Hence since we have a continuous function over a bounded interval, one can conclude via Lemma \ref{Lem:Swap derivative and integral} and automatic differentiation,
	\begin{align*}
	\partial_{x}^{i} \hat{v}(t,x)
	& 
	=
	\partial_{x}^{i} \bE \Big[ \eta(T_{1}, \bar{X}_{T_{1}}) \big| \cF_{t} \Big]
	 =
	\bE \left[ \eta(T_{1}, \bar{X}_{T_{1}}) \cW_{1}^{i} \big| \cF_{t} \right] \notag
	=
	\bE \left[ \psi^{t,x} \cW_{1}^{i} \big| \cF_{t} \right] \, .
	\end{align*}
	Technically we have again used the Tower property to remove the final conditional expectation which requires integrability. We now show this is valid and due to the form of $\psi$ we split into two terms,
	\begin{align*}
	&  \bE \left[ \psi^{t,x} \cW_{1}^{i} \big| \cF_{t} \right]
	 =
	\bE \left[ 
	\1_{\{N_{T} =0 \}} \psi^{t,x} \cW_{1}^{i}
	+
	\1_{\{N_{T} \ge 1\}} \psi^{t,x} \cW_{1}^{i}
	\big| \cF_{t} \right] \, .
	\end{align*}
	One can automatically see that if $N_{T} \ge 1$ then $\psi=\tilde{\psi}$, for the case $N_{T} =0$, we need to show equivalence between $\psi$ and the corresponding $\Phi$. Firstly let us show,
	\begin{align*}
	\bE \left[ 
	\1_{\{N_{T} =0 \}} \psi^{t,x} \cW_{1}^{1}
	\big| \cF_{t} \right]
	=
	\bE \left[ 
	\1_{\{N_{T} =0 \}} \Phi_{1}^{t,x} \cW_{1}^{1} 
	\big|\cF_{t} \right] \, .
	\end{align*}
	Expanding out $\Phi_{1}$ we obtain,
	\begin{align*}
	&\bE \left[ 
	\1_{\{N_{T} =0 \}} \Phi_{1}^{t,x} \cW_{1}^{1}
	\big| \cF_{t} \right] 
	=
	\bE \left[ 
	\1_{\{N_{T} =0 \}} \frac{g(\bar{X}_{T_{N_{T}+1}})- g(\hat{X}_{T_{N_{T}+1}})}{2\overline{F}(\Delta T_{N_{T}+1})} \cW_{1}^{1} 
	\Big| \cF_{t} \right] 
	\, .
	\end{align*}
	Using that $W$ and $-W$ have the same distribution and $\cW^{1}$ is an odd function of the Brownian increment $\Delta W$ (see \eqref{Eq:First and Second Malliavin Weights}) we obtain,
	\begin{align*}
	&\bE \left[ 
	\1_{\{N_{T} =0 \}} \Phi_{1}^{t,x} \cW_{1}^{1}
	\big| \cF_{t} \right] 
	=
	\bE \left[ 
	2\1_{\{N_{T} =0 \}} \frac{g(\bar{X}_{T_{N_{T}+1}})}{2\overline{F}(\Delta T_{N_{T}+1})} \cW_{1}^{1} 
	\Big| \cF_{t} \right] 
	\, ,
	\end{align*}
	which shows the required result. Equivalently, we now show the equality
	\begin{align*}
	\bE \left[ 
	\1_{\{N_{T} =0 \}} \psi^{t,x} \cW_{1}^{2} 
	\big| \cF_{t} \right]
	=
	\bE \left[ 
	\1_{\{N_{T} =0 \}} \Phi_{2}^{t,x} \cW_{1}^{2}
	\big|\cF_{t} \right] \, .
	\end{align*}
	By a similar argument to above,
	\begin{align*}
	&\bE \left[ 
	\1_{\{N_{T} =0 \}} \Phi_{2}^{t,x} \cW_{1}^{2}
	\big| \cF_{t} \right] 
	\\
	& =
	\bE \left[ 
	\1_{\{N_{T} =0 \}} \frac{g(\bar{X}_{T_{N_{T}+1}})+ g(\hat{X}_{T_{N_{T}+1}})
		-2 g(\bar{X}_{T_{N_{T}}} + b(T_{N_{T}}) \Delta T_{N_{T}+1})
	}{2\overline{F}(\Delta T_{N_{T}+1})} \cW_{1}^{2}
	\Big| \cF_{t} \right] 
	\, .
	\end{align*}
	By the fact that $g(\bar{X}_{T_{N_{T}}} + b(T_{N_{T}}) \Delta T_{N_{T}+1})$ is $\overline{\cF}_{N_{T}}$-adapted, and the weight has zero expectation we can remove this term from the expectation. Again, since $W$ and $-W$ have the same distribution, and $\cW^{2}$ is even we obtain,
	\begin{align*}
	&\bE \left[ 
	\1_{\{N_{T} =0 \}} \Phi_{2}^{t,x} \cW_{1}^{2}
	\big| \cF_{t} \right] 
	=
	\bE \left[ 
	2\1_{\{N_{T} =0 \}} \frac{g(\bar{X}_{T_{N_{T}+1}})}{2\overline{F}(\Delta T_{N_{T}+1})} \cW_{1}^{2}
	\Big| \cF_{t} \right] 
	\, ,
	\end{align*}
	again, this yields the required result. Thus the spatial derivatives of $\hat{v}$ satisfy,
	\begin{align*}
	\partial_{x}^{i} \hat{v}(t,x)
	& =
	\bE \left[ 
	\1_{\{N_{T} =0 \}} \Phi_{i}^{t,x} \cW_{1}^{i}
	+
	\1_{\{N_{T} \ge 1\}} \tilde{\psi}^{t,x} \cW_{1}^{i}
	\big| \cF_{t} \right].
	\end{align*}
	Uniform integrability of $\tilde{\psi} \cW^{i}$ and $\Phi_{i} \cW^{i}$ then implies $\partial_x^{i} \hat{v}(t,x)$ is a continuous function and one can use this integrability to also conclude $\partial_{x}^{i} \hat{v}(t,x) = 	\bE \left[ \psi^{t,x} \cW_{1}^{i} \big| \cF_{t} \right]$. Thus existence of the first and second spatial derivatives are assured.

	\emph{Step 4: Representations match.}
	Introducing the following notation, $N_{T}(s):=N_{T}-N_{s}$, i.e. the number of regime switches that occur between time $s$ and $T$, with the obvious relation $N_{T}(t)=N_{T}$.
	
	To show that the two representations are the same, we need to consider the terms $\partial_{x}^{i} \hat{v}(T_{1}, \bar{X}_{T_{1}})$ for $t \le T_{1} < T$. One has that,
	\begin{align*}
	\hat{v}(T_{1}, \bar{X}_{T_{1}}) = \bE[ \psi^{T_{1}, \bar{X}_{T_{1}}} | \cF_{T_{1}}] \, .
	\end{align*}
	To apply derivatives we again introduce the function $\eta(T_{2}, \bar{X}_{T_{2}}) =\bE[ \psi^{T_{1}, \bar{X}_{T_{1}}} | \cF_{T_{2}}]$ and then Lemma \ref{Lem:Swap derivative and integral} and Malliavin automatic differentiation implies,
	\begin{align*}
	\partial_{x}^{i} \hat{v}(T_{1}, \bar{X}_{T_{1}}) 
	=
	 \bE[ \psi^{T_{1}, \bar{X}_{T_{1}}} \cW_{2}^{i} | \cF_{T_{1}}]
	 \quad
	 \bP\text{-a.s.}
	\end{align*}
	Using the same arguments as before we can rewrite this as,
	\begin{align*}
	&\partial_{x}^{i} \hat{v}(T_{1}, \bar{X}_{T_{1}})
	 =
	\bE \left[ \1_{\{N_{T}(T_{1}) =0 \}} \Phi_{i}^{T_{1}, \bar{X}_{T_{1}}} \cW_{2}^{i}
	+
	\1_{\{N_{T}(T_{1}) \ge 1\}} \tilde{\psi}^{T_{1}, \bar{X}_{T_{1}}} \cW_{2}^{i} \Big| \cF_{T_{1}} \right] 
	\quad
	\bP\text{-a.s.}
	\end{align*}
	One then recognises the internal conditional expectations in \eqref{Eq:At Least One Branching} as the derivatives of $\hat{v}$ starting at time $(T_{1}, \bar{X}_{T_{1}})$. Thus, by integrability, \eqref{Eq:At Least One Branching} can be simply written as,
	\begin{align*}
	\bE \left[
	\1_{\{ N_{T} \ge 1 \}}
	\frac{1}{f(\Delta T_{1})}
	\Bigg(
	\Delta b_{1}
	\partial_{x}
	\hat{v}(T_{1},\bar{X}_{T_{1}})
	-
	\frac{1}{2} \sigma(\theta_{0},T_{1})^{2}
	\partial_{xx}
	\hat{v}(T_{1},\bar{X}_{T_{1}})
	\Bigg)
	~
	\Bigg | ~
	\cF_{t}
	\right] \, .
	\end{align*}
	This leads us to the following nonlinear relation for $\hat{v}$,
	\begin{align*}
	\hat{v}(t,x)= \bE \left[
	\frac{g(\bar{X}_{T_{1}})}{\overline{F}(\Delta T_{1})} \1_{\{ N_{T}=0 \}}
	+
	\1_{\{ N_{T} \ge 1 \}}
	\frac{
		\Delta b_{1}
		\partial_{x}
		\hat{v}(T_{1},\bar{X}_{T_{1}})
		-
		\frac{1}{2} \sigma(\theta_{0},T_{1})^{2}
		\partial_{xx}
		\hat{v}(T_{1},\bar{X}_{T_{1}})
	}{f(\Delta T_{1})}
	~
	\Bigg | ~
	\cF_{t}
	\right] \, .
	\end{align*}
	Since this representation and \eqref{Eq:Forward Representation} are equal we have $v(t,x) = \hat{v}(t,x)$ hence our representation solves the PDE.
\end{proof}


\subsection{Verifying the integrability assumptions}
Theorem \ref{theo:StochasticRepresentation} relied on various integrability assumptions and our final result is to show that these assumptions hold.

\begin{theorem}
	\label{Thm:Integrability Hold}
	Let Assumptions \ref{Assume:Main SDE Assumptions}, \ref{Assume:Bounded g} and \ref{Assume:Finite Variance} hold. Then the integrability conditions in Theorem \ref{theo:StochasticRepresentation} hold.
\end{theorem}

\begin{proof}
	We start by showing the uniform integrability conditions, recall that for uniform integrability to hold it is sufficient to show the stochastic process is in $L^{p}$ for $p>1$ (see \cite{Williams1991}*{Chapter 13} for results on uniform integrability).
	
	Firstly, by Proposition \ref{Prop:Representation}, one can conclude that $\psi^{t,x} \in L^{2}$, thus we have the required uniform integrability. Let us now consider $\tilde{\psi}^{t,x} \cW_{1}^{1}$ and $\tilde{\psi}^{t,x} \cW_{1}^{2}$. Due to both quantities having very similar forms we consider $\tilde{\psi}^{t,x} \cW_{1}^{i}$ for $i \in \{1, 2\}$, hence we want to show,
	\begin{align*}
	\bE[ |\tilde{\psi}^{t,x} \cW_{1}^{i}|^{p} | \cF_{t}] < \infty,
	\quad
	\text{for some }
	p >1.
	\end{align*}
	We show this by borrowing many of the arguments in the proof of Proposition \ref{Prop:Representation}, hence we take $p=2$. Using the representation for $\tilde{\psi}^{t,x}$ and taking common factors we obtain,
	\begin{align*}
	\bE[ |\tilde{\psi}^{t,x} \cW_{1}^{i}|^{2} | \cF_{t}] 
	\le
	& \bE \Bigg[
	 \Bigg(
	 \frac{\Delta g_{T_{N_{T}+1}} - \Delta \hat{g}_{T_{N_{T}+1}}}{2\overline{F}(\Delta T_{N_{T}+1})}
	 \frac{\Delta b_{N_{T}} \cW_{N_{T}+1}^{1}}{f(\Delta T_{N_{T}})} 
	 \Bigg)^{2}
	 \prod_{k=2}^{N_{T}}
	 P_{k}^{2}
	 \left(\cW_{1}^{i}\right)^{2}
	 \Big | ~ \cF_{t}
	 \Bigg]
	 \\
	 &
	 +
	 \bE \Bigg[
	 \Bigg(
	 \frac{\Delta g_{T_{N_{T}+1}} + \Delta \hat{g}_{T_{N_{T}+1}}}{2\overline{F}(\Delta T_{N_{T}+1})}
	 \frac{\frac{1}{2} \sigma(\theta_{N_{T}-1},T_{N_{T}})^{2} \cW_{N_{T}+1}^{2}}{f(\Delta T_{N_{T}})} 
	 \Bigg)^{2}
	 \prod_{k=2}^{N_{T}}
	 P_{k}^{2}
	 \left(\cW_{1}^{i}\right)^{2}
	 \Big | ~ \cF_{t}
	 \Bigg] \, .
	\end{align*}
	We now use the same techniques from the proof of Proposition \ref{Prop:Representation}, firstly, we can condition on $N_{T}= \ell$ and multiply by the corresponding probability. Then by conditioning on $\overline{\cF}_{N_{T}}$ (see proof of Proposition \ref{Prop:Representation}) we obtain the following,
	\begin{align*}
	\bE \Bigg[
	\Bigg(
	\frac{\Delta g_{T_{N_{T}+1}} - \Delta \hat{g}_{T_{N_{T}+1}}}{2\overline{F}(\Delta T_{N_{T}+1})}
	\frac{\Delta b_{N_{T}} \cW_{N_{T}+1}^{1}}{f(\Delta T_{N_{T}})} 
	\Bigg)^{2}
	\Big | ~ \overline{\cF}_{N_{T}}, ~ N_{T}= \ell
	\Bigg]
	\le
	C
	\frac{\Delta b_{N_{T}}^{2}}{f(\Delta T_{N_{T}})^{2}} \, ,
	\end{align*}
	and
	\begin{align*}
	& \bE \Bigg[
	\Bigg(
	\frac{\Delta g_{T_{N_{T}+1}} + \Delta \hat{g}_{T_{N_{T}+1}}}{2\overline{F}(\Delta T_{N_{T}+1})}
	\frac{\frac{1}{2} \sigma(\theta_{N_{T}-1},T_{N_{T}})^{2} \cW_{N_{T}+1}^{2}}{f(\Delta T_{N_{T}})} 
	\Bigg)^{2}
	\Big | ~ \overline{\cF}_{N_{T}}, ~ N_{T}= \ell
	\Bigg]
	\\
	&
	\qquad
	\le
	\frac{
		C
	}{
	f(\Delta T_{N_{T}})^{2}}
\frac{\sigma(\theta_{N_{T}-1},T_{N_{T}})^{4}}{\sigma(\theta_{N_{T}},T_{N_{T}+1})^{2}}
\frac{\min \left[
	1
	,
	\sigma(\theta_{N_{T}},T_{N_{T}+1})^{2} \Delta T_{N_{T}+1}
	\right]}{\Delta T_{N_{T}+1}}
	 \, .
	\end{align*}
	We now use these bounds to bound $\tilde{\psi}\cW$. Concentrating on the $\Delta b_{N_{T}}$ term, we follow the finite variance proof and condition out $\Delta b_{N_{T}}^{2} P_{N_{T}}^{2}$, then use \eqref{Eq:Bound b P term}, namely,
	\begin{align*}
	 & \bE \Bigg[
	\Bigg(
	\frac{\Delta g_{T_{N_{T}+1}} - \Delta \hat{g}_{T_{N_{T}+1}}}{2\overline{F}(\Delta T_{N_{T}+1})}
	\frac{\Delta b_{N_{T}} \cW_{N_{T}+1}^{1}}{f(\Delta T_{N_{T}})} 
	\Bigg)^{2}
	\prod_{k=2}^{N_{T}}
	P_{k}^{2}
	\left(\cW_{1}^{i}\right)^{2}
	\Big | ~ \cF_{t}, ~ N_{T}= \ell
	\Bigg]
	\\
	&
	\quad
	\le
	\bE \Bigg[
	\frac{C}{f(\Delta T_{N_{T}})^{2}}
	\frac{1}{f(\Delta T_{N_{T}-1})^{2}} 
	\frac{\sigma(\theta_{N_{T}-2},T_{N_{T}-1})^{4}}{\sigma(\theta_{N_{T}-1},T_{N_{T}})^{4}}
	\prod_{k=2}^{N_{T}-1}
	P_{k}^{2}
	\left(\cW_{1}^{i}\right)^{2}
	\Big | ~ \cF_{t}, ~ N_{T}= \ell
	\Bigg]
	\, .
	\end{align*}
	By continuing to follow the argument we can bound the above quantity by,
	\begin{align}
	\label{Eq:b Integrability bound}
	\bE \Bigg[
	\frac{C^{N_{T}}}{f(\Delta T_{N_{T}})^{2}}
	\frac{\sigma(\theta_{0},T_{1})^{4}}{\sigma(\theta_{N_{T}-1},T_{N_{T}})^{4}}
	\prod_{k=2}^{N_{T}-1}
	\frac{1}{f(\Delta T_{k})^{2} \Delta T_{k}^{2}}
	\frac{1}{f(\Delta T_{1})^{2}} 
	\left(\cW_{1}^{i}\right)^{2}
	\Big | ~ \cF_{t}, ~ N_{T}= \ell
	\Bigg]
	\, .
	\end{align}
	Since $\sigma_{0}>0$ is constant it is clear that,
	\begin{align*}
	\bE[\left(\cW_{1}^{1}\right)^{2}| \overline{\cF}_{0}]
	\le
	C \bE[\left(\cW_{1}^{2}\right)^{2}| \overline{\cF}_{0}]
	\le
	C \frac{1}{\Delta T_{1}^{2}} \, .
	\end{align*}
	Hence we can bound \eqref{Eq:b Integrability bound},
	\begin{align*}
	\bE \Bigg[
	\frac{C^{N_{T}}}{f(\Delta T_{N_{T}})^{2}}
	\frac{\sigma(\theta_{0},T_{1})^{4}}{\sigma(\theta_{N_{T}-1},T_{N_{T}})^{4}}
	\prod_{k=1}^{N_{T}-1}
	\frac{1}{f(\Delta T_{k})^{2} \Delta T_{k}^{2}}
	\Big | ~ \cF_{t}, ~ N_{T}= \ell
	\Bigg]
	\le
	\bE \Big[
	C^{N_{T}}
	\Big | ~ \cF_{t}, ~ N_{T}= \ell
	\Big]
	\, ,
	\end{align*}
	where the inequality follows from our assumptions on $f$ and $\sigma$.
	
	Using this argument to deal with the extra Malliavin weight and the arguments in Proposition \ref{Prop:Representation}, we also obtain,
	\begin{align*}
	& \bE \Bigg[
	\Bigg(
	\frac{\Delta g_{T_{N_{T}+1}} + \Delta \hat{g}_{T_{N_{T}+1}}}{2\overline{F}(\Delta T_{N_{T}+1})}
	\frac{\frac{1}{2} \sigma(\theta_{N_{T}-1},T_{N_{T}})^{2} \cW_{N_{T}+1}^{2}}{f(\Delta T_{N_{T}})} 
	\Bigg)^{2}
	\prod_{k=2}^{N_{T}}
	P_{k}^{2}
	\left(\cW_{1}^{i}\right)^{2}
	\Big | ~ \cF_{t}
	\Bigg] 
	 \\
	 &
	 \quad
	 \le
	 \bE \left[ C^{N_{T}}
	 \frac{\sigma(\theta_{N_{T}},T_{N_{T}+1})^{\nu} \Delta T_{N_{T}+1}^{\nu /2}}
	 {\sigma(\theta_{N_{T}},T_{N_{T}+1})^{2}} \Delta T_{N_{T}+1}^{-1}
	 \prod_{k=1}^{N_{T}}
	 \Delta T_{k}^{-2 \kappa}
	 \Bigg | N_{T}=\ell \right] ,
	 \quad \text{for } \nu \in [0,2]
	\, .
	\end{align*}
The finiteness of these bounds follows directly from Proposition \ref{Prop:Representation}.
	
For the $f(\Delta T_{1})^{-1}\Delta b_{1} \tilde{\psi}^{T_{1},\bar{X}_{T_{1}}} \cW_{2}^{1} $ and $f(\Delta T_{1})^{-1}\sigma(\theta_{0}, T_{1})^{2}  \tilde{\psi}^{T_{1},\bar{X}_{T_{1}}} \cW_{2}^{2} $ terms, these follow automatically from Proposition \ref{Prop:Representation}.

	For uniform integrability of $\Phi_{1} \cW^{1}$, take $p=2$ as above. Then use Cauchy-Schwarz and the Lipschitz property of $g$, which yields $| \Delta g_{T_{N_{T}+1}} - \Delta \hat{g}_{T_{N_{T}+1}}| \le C \sigma(\theta_{N_{T}}, T_{N_{T}+1}) | \Delta W_{T_{N_{T}+1}}|$. One notes that the $\sigma$ and $\Delta T$ terms cancel and hence finite.
	
	Similarly, for $\Phi_{2} \cW^{2}$, again take $p=2$ and use Cauchy-Schwarz along with \eqref{Eq:Ito Bound on g}. Again all terms cancel which implies this is also finite and hence uniformly integrable.
	
	The final integrability results we require are all $\bP$-a.s. results. We have $\psi^{T_{1}, \bar{X}_{T_{1}}}$, $\Delta b_{2} \tilde{\psi}^{T_{2},\bar{X}_{T_{2}}} \cW_{3}^{1}$ and
	$\sigma(\theta_{1}, T_{2})^{2}  \tilde{\psi}^{T_{2},\bar{X}_{T_{2}}} \cW_{3}^{2}$ are $\bP$-a.s. uniformly integrable, and $\tilde{\psi}^{T_{1},\bar{X}_{T_{1}}} \cW_{2}^{1}$ and $\tilde{\psi}^{T_{1},\bar{X}_{T_{1}}} \cW_{2}^{2}$ are $\bP$-a.s. integrable. However, these follow from the arguments above along with the fact that $t<T_{1} < T_{2}$ $\bP$-a.s. hence $\sigma(\theta_{1}, T_{2}) < \infty$ $\bP$-a.s.
	 Hence we have shown all the required integrability conditions to use Theorem \ref{theo:StochasticRepresentation}. 	
\end{proof}

The proof of Theorem \ref{Thm:Representation Solves the PDE} follows in a straightforward way by combining these results.
\begin{proof}[Proof of Theorem \ref{Thm:Representation Solves the PDE}]
	By letting Assumptions \ref{Assume:Main SDE Assumptions}, \ref{Assume:Bounded g} and \ref{Assume:Finite Variance} hold, then Theorems \ref{theo:StochasticRepresentation} and \ref{Thm:Integrability Hold} imply that our estimator $\tilde{v}$ given in \eqref{Eq:Variance Reduced Representation} solves the PDE \eqref{Eq:High Dimension Transport PDE}.
	
	Moreover, Proposition \ref{Prop:Representation}, implies that $\psi$ is square integrable and hence of finite variance.
	
\end{proof}


\section{Towards the general case and future work}
\label{Sec:Toward the general case}
The methodology presented in this work can be extended to accommodate PDEs of the form,
\begin{equation}
\label{Eq:More General Transport PDE}
\begin{cases}
\partial_{t} v(t,x) + b(t)\cdot Dv(t,x) + h(t,x) =0  \quad\text{for all }
(t,x) \in [0,T)\times \bR^{d} \, ,
\\
v(T,x)=g(x) \, ,
\end{cases}
\end{equation}
where $h$ is a nice function and we still have $v \in C_{b}^{1,3}$. As in the case of standard branching representations one introduces a further probability measure $\bP_{B}$ on the space $\{ 0, 1 \}$, where $0$ signifies the case the particles dies (this can be thought of as a $v^{0}$ term) at position $(T_{k}, \bar{X}_{T_{k}})$ and we evaluate $h$ at this position.

\subsection{Allowing the drift to have a spatial dependence}
\label{Sec:b with space component}
Throughout this chapter we have made the assumption that the drift $b$ does not depend on space. The main reason for this is to ensure finite variance. One can consider replacing Assumption \ref{Assume:Main SDE Assumptions}, with $b: [0,T] \times \bR \rightarrow \bR$, satisfying $1/2$-H\"{o}lder in time, Lipschitz in space and uniformly bounded and most of the arguments presented still hold. The bound that changes and makes the arguments more difficult is \eqref{Eq:Full P Bound}, to see this let us observe how $\Delta b$ is bounded under these new assumptions,
\begin{align*}
\bE\big[\Delta b_{k}^{4} & | \overline{\cF}_{k-1}, N_{T}= \ell\big] 
\\ &
\le
C \bE\Big[(b(T_{k},\bar{X}_{T_{k}})-b(T_{k},\bar{X}_{T_{k-1}}) )^{4} + (b(T_{k},\bar{X}_{T_{k-1}})-b(T_{k-1},\bar{X}_{T_{k-1}}) )^{4} | \overline{\cF}_{k-1}, N_{T}= \ell\Big] \, .
\end{align*}
For the second term we can use $1/2$-H\"{o}lder continuity in time of $b$, for the first term we can Lipschitz continuity in space to obtain,
\begin{align*} 
\bE[(b(T_{k},\bar{X}_{T_{k}})-b(T_{k},\bar{X}_{T_{k-1}}) )^{4} | \overline{\cF}_{k-1}, N_{T}= \ell] 
\le &
C\bE[(\bar{X}_{T_{k}} - \bar{X}_{T_{k-1}} )^{4} | \overline{\cF}_{k-1}, N_{T}= \ell] 
\\
\le &
C \bE[( \Delta T_{k} + \sigma(\theta_{k-1},T_{k}) \Delta W_{T_{k}} )^{4} | \overline{\cF}_{k-1}, N_{T}= \ell] 
\\
\le &
C \sigma(\theta_{k-1},T_{k})^4 \Delta T_{k}^{2} 
\, .
\end{align*}
Since $\sigma$ is bounded from below we can conclude,
\begin{align*}
\bE[\Delta b_{k}^{4} | \overline{\cF}_{k-1}, N_{T}= \ell] 
\le 
C \sigma(\theta_{k-1},T_{k})^4 \Delta T_{k}^{2}  \, .
\end{align*}
It is also straightforward to see the same bound applies if we take $b$ Lipschitz in time. The bounds on $M$ and $V$ still have the form
\begin{align*}
\bE[M_{k+1}^{4}| \overline{\cF}_{k}, N_{T}= \ell]
& \le
C
\frac{\Delta b_{k}^{4}}{\Delta T_{k+1}^{2}} \frac{1}{\sigma(\theta_{k},T_{k+1})^{4}}
\, , 
\\
\bE[V_{k+1}^{4}| \overline{\cF}_{k}, N_{T}= \ell]
& \le
C \frac{\sigma(\theta_{k-1},T_{k})^{8}}{\sigma(\theta_{k},T_{k+1})^{8}}
\frac{1}{\Delta T_{k+1}^{4}}
\, ,
\end{align*}
although one should note that we cannot use the $\Delta b$ bound above in the $M$ term since they are w.r.t. different conditional expectations. That being said though one can still observe where a problem arises by considering,
\begin{align*}
& \bE\big[ \bE[P_{k+1}^{4}| \overline{\cF}_{k}, N_{T}= \ell] \big| \overline{\cF}_{k-1}, N_{T}= \ell \big]
\\
& \le 
C \bE \Big[ \frac{1}{f(\Delta T_{k})^{4}} 
\Big(
\frac{\Delta T_{k}^{2} }{\Delta T_{k+1}^{2}} 
\frac{\sigma(\theta_{k-1},T_{k})^4 }{\sigma(\theta_{k},T_{k+1})^{4}}
+
\frac{\sigma(\theta_{k-1},T_{k})^{8}}{\sigma(\theta_{k},T_{k+1})^{8}}
\frac{1}{\Delta T_{k+1}^{4}}
\Big)
\Big| \overline{\cF}_{k-1}, N_{T}= \ell \Big].
\end{align*} 
Whereas in the proof we can bound \eqref{Eq:Full P Bound} by the term arising from the $V$ (i.e. the $V$ bound dominates the $M$ bound), that is not the case here. To see this take $n=-1$ for the coefficient in the $\sigma$, we then obtain,
\begin{align*}
C \bE \Big[ \frac{1}{f(\Delta T_{k})^{4}} 
\frac{\Delta T_{k}^{6} }{\Delta T_{k+1}^{2}} 
\Big(
1
+
\frac{\Delta T_{k}^{2}}{\Delta T_{k+1}^{2}}
\Big)
\Big| \overline{\cF}_{k-1}, N_{T}= \ell \Big].
\end{align*} 
Therefore the $1$ (term arising from the $M$) is larger if $\Delta T_{k} < \Delta T_{k+1}$, hence we cannot dominate in the same way. As it turns out this a not a problem for obtaining \eqref{Eq:Bound for b term}, however, it does become an issue for obtaining \eqref{Eq:Full bound with nu}. This appears because \eqref{Eq:Full bound with nu} relies on a cancelling argument, while this extra term changes the original bound from,
\begin{align*}
\bE \Big[ \frac{1}{f(\Delta T_{k})^{4}} 
\frac{\Delta T_{k}^{8} }{\Delta T_{k+1}^{4}} 
\Big| \overline{\cF}_{k-1}, N_{T}= \ell \Big]
\quad
\text{to}
\quad
\bE \Big[ \frac{1}{f(\Delta T_{k})^{4}} 
\frac{\Delta T_{k}^{8} }{\Delta T_{k+1}^{4}\Delta T_{k}^{2}} 
\Big| \overline{\cF}_{k-1}, N_{T}= \ell \Big].
\end{align*}
This extra $\Delta T_{k}$ dependency makes the bound far weaker and consequently proving finite variance becomes more difficult. Of course the new bound we have obtained is not sharp, for example in the case $\Delta T_{k} \ge \Delta T_{k+1}$ we can return the original bound. 

If we wish to argue the proof in a similar way one must either look to obtain a stronger bound on $\Delta b$ (this is essentially why $b$ in Assumption \ref{Assume:Main SDE Assumptions} worked), or one can find a way to make the $V$ term dominate without increasing its size so much to break the remainder of the proof. For example, an interesting route to explore is to add an event probability distribution to the $M$ and $V$ term (similar to other branching diffusion algorithms) applying a judicious choice of probability distribution may give us the means to bound the $M$ term by $V$ again. 

There are of course many different approaches one can take to solve this problem and as described, the remaining arguments in Theorems \ref{theo:StochasticRepresentation} and \ref{Thm:Integrability Hold} follow with a more general $b$. But proving finite variance of this representation remains an open question.

\subsection{Fully nonlinear first order case}
Of course the true end goal of this work is to handle nonlinearities, for example, Burger's type $v Dv$, which arise in many applications and for which numerical methods like characteristics cannot apply. Therefore future work will be on addressing explicit conditions under which this method provides solutions to transport PDEs of the form,
\begin{equation*}
\begin{cases}
\partial_{t} v(t,x) + b(t,x) \cdot Dv(t,x) =f(t,x,v,Dv) \, ,
\\
v(T,x)=g(x) \, ,
\end{cases}
\end{equation*}
where $f$ is polynomial in $v$ and $Dv$. 

Handling such general first order PDEs will require additional arguments to what we have presented here. However, ideas from the case $b(t,x)$ along with the (purely numerical) technique presented in \cite{Warin2017} may yield the necessary tools to overcome such equations.

\begin{remark}[Requirement for Smooth Solutions]
	In theory this technique should be able to extend to the general, fully nonlinear case, one will still require a sufficiently smooth classical solution to the underlying PDE. The reason for this is due to the fact we assign a representation to $\partial_{xx}v$, thus we automatically require existence of this quantity.
	
	This implies that if we argue that the representation solves the PDE via viscosity solutions then we in fact show a classical solution. Of course this implies the method is not suitable for PDEs with ``shocks''.  
\end{remark}


\section{Examples}
\label{sec:examples}
We show the potential of this method on two examples to compare this technique against the standard perturbation technique. 
The first example is a simple linear PDE which satisfies all of our assumptions and hence is only an example to show that our algorithm converges to the true, while the perturbation converges to a different value. The second is a nonlinear first order PDE, this is the more interesting case and we still observe our method giving reasonable results.  

\subsection{Simple First Order PDE}
\begin{figure}[!ht] 
	\centering
	\includegraphics[width=15cm]{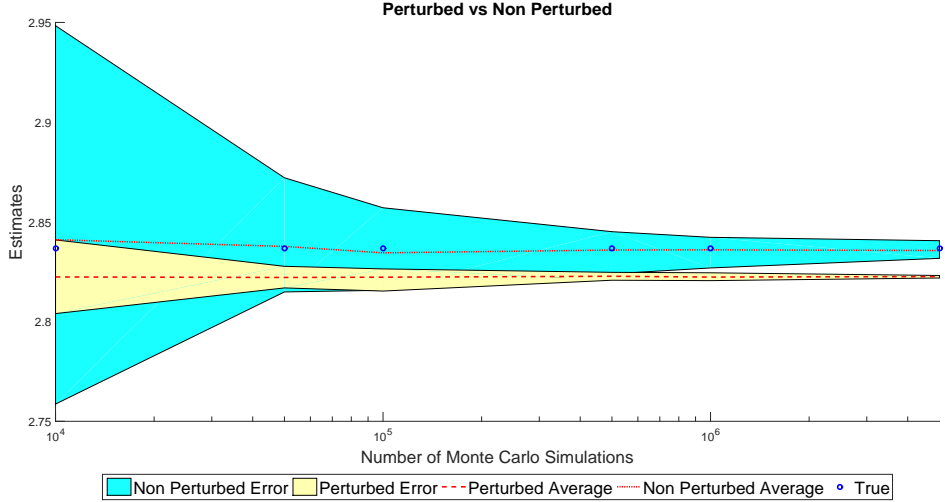}
	\vspace{-1.0cm}
	\caption{Shows the error and estimates of the solution as a function of the number of Monte Carlo simulations. The error corresponds to the approximate 90\% confidence interval.}
	\label{Fig:Basic PDE example for Bias}
\end{figure}

Let us consider the following linear PDE,
\begin{equation}
\label{Eq:Simple Example PDE}
\begin{cases}
\partial_{t} v(t,x) + \partial_{x}v(t,x) =0  \quad\text{for all }
(t,x) \in [0,1)\times \bR \, ,
\\
v(1,x)= 10\cos(x-1-5) \, .
\end{cases}
\end{equation}
It is then clear to see that $v(t,x)=10\cos(x-t-5)$ satisfies this PDE. Although such a PDE is easy to solve it serves as a good example to show the issue using a perturbation. We want to solve this PDE at the point $(0,10)$, where the true solution is $\approx 2.84$. By considering the case where we perturb by $\sigma=0.1$, and then estimate the expectation using varying amounts of Monte Carlo simulations, see Figure \ref{Fig:Basic PDE example for Bias}. To get a handle on the variance (error) we ran the simulation 50 times, plotted the average and the approximate 90\% confidence interval. That is we view the largest and smallest value as a proxy for convergence of the algorithm. For the unbiased algorithm we also took, $n=-1$ and for the Gamma parameters $\kappa=1/2$ and $\eta=2$.

What is clear from Figure \ref{Fig:Basic PDE example for Bias} is, as the number of Monte Carlo simulations increase, both algorithms are converging. However the perturbed case stays at a constant level away from the true, which implies that the estimate is biased (as was expected). Therefore no amount of Monte Carlo simulations will yield the true solution. For the unbiased algorithm, although having a higher variance, we see that the average hovers around the true value and moreover we observe convergence towards this point.

Hence the stochastic representation we derive indeed yields the true solution of the PDE, what is more fascinating and important about this result though is $\sigma$ is not tending to zero, in fact we can bound it from below, this is the key step when it comes to more complex PDEs.

Moreover, this calculation was carried out using a basic Monte Carlo algorithm, one could look to more sophisticated techniques as appearing in \cite{DoumbiaEtAl2017} where the authors apply particle methods for an improved convergence.

\subsection{Nonlinear PDE}
Let us now generalise to the nonlinear setting and consider the following PDE,
\begin{equation}
\label{Eq:Nonlinear PDE}
\begin{cases}
\partial_{t} v(t,x) + \partial_{x}v(t,x) +\frac{1}{10} \big( (\partial_{x}v(t,x))^2 + v(t,x)^2 -1 \big) =0  \quad\text{for all }
(t,x) \in [0,1)\times \bR \, ,
\\
v(1,x)= \cos(1-x) \, .
\end{cases}
\end{equation}
We have taken this PDE since it is simple to observe that $v(t,x)=\cos(t-x)$ is the solution. It also is nice enough that one would expect our unbiased algorithm and the perturbation algorithm to work reasonably well. We want to solve this at the point $(0,1)$.

$\triangleright$ \emph{Convergence issue for the perturbation algorithm} One can note that, applying the perturbation technique implies that the resulting PDE is a second order semilinear PDE, and hence the corresponding branching algorithm is given in \cite{LabordereEtAl2016}. This creates a problem for the convergence of the algorithm, Assumption 3.10 and Theorem 3.12 of \cite{LabordereEtAl2016} give minimum bounds on the relative size of the drift to the diffusion, even for \eqref{Eq:Nonlinear PDE} which has a extremely nice solution, we observe that the algorithm fails to converge for $\sigma_{0}=0.5$ and has a large variance for $\sigma_{0}$ smaller than 1. Needless to say this is not a desirable property for the algorithm to have; perturbation can only work as a method if the perturbation is small and here we observe that there is a lower bound on the size of the perturbation and hence the bias of the estimator. Furthermore, as it turns out, there is no such problem with our unbiased algorithm and one can observe convergence for $\sigma_{0} < 0.5$.

With the above in mind, in order to make the two algorithms comparable we set the perturbed algorithm as $\sigma_{0}=1$, but the remaining parameters are as above. Because the variance here is larger than the linear PDE we consider 100 realisations for each Monte Carlo level then take the approximate 80\% confidence intervals and the average is then based on these 80 realisations. Furthermore, because we are dealing with nonlinear terms we have a more complex representation and need to establish a probability distribution for the type of event i.e. $v^2$, $(\partial_{x}v)^2$ etc. This is well understood in the case of the perturbation algorithm (see \cite{LabordereEtAl2016}), however, the variance of our unbiased algorithm seems to be highly dependent on how one chooses this probability distribution.

\begin{figure}[!ht] 
	\centering
	\includegraphics[width=15cm]{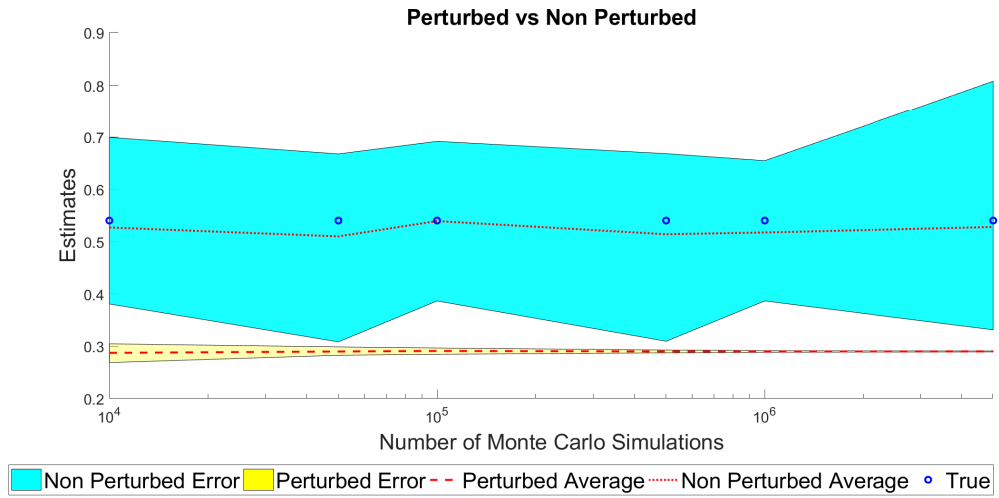}
	\vspace{-1.0cm}
	\caption{Shows the error and estimates of the solution as a function of the number of Monte Carlo simulations. The error corresponds to the approximate 80\% confidence interval.}
	\label{Fig:Nonlinear PDE example}
\end{figure}

Figure \ref{Fig:Nonlinear PDE example} shows that yet again our unbiased algorithm provides a correction for the second order term. While the perturbation algorithm converges to a different value. However, it is clear that the variance in our algorithm is much higher. One of the reasons for this is because of the uncertainty in what events will be used for each realisation. Namely, for the linear PDE case, there was no probability distribution over events and this allowed us to bound the variance. In this more general case, more work would have to be done in order to bound the variance, and from our numerical example the choice of probability distribution has a role to play here.

\section*{Conclusions and Outlook}
We have demonstrated a stochastic algorithm capable of dealing with first order PDEs, where originally such PDEs seemed beyond the reach of stochastic methods without approximation. This has potentially large implications for numerics of such PDEs since stochastic algorithms can easily be parallelised and scale favourable with dimension as argued in \cite{BernalEtAl2017}.

Due to the added difficulty in considering more general transport PDEs we have taken a simple case here. As a consequence we have left some open problems to be addressed, namely.
\begin{enumerate}
	\item Finite variance estimator when the drift component also depends on space.
	
	\item Dropping the assumption on the initial PDE having a classical solution.
	
	\item Extending to the case of nonlinear terms in both the solution of the PDE and its first spatial derivative.
\end{enumerate}

Our hope is that with the continued research and innovation into branching diffusions that such results will be within reach.

\bigskip

\noindent\textit{Acknowledgements.}
The authors would like to thank Nizar Touzi (\'Ecole Polytechnique Paris) and Christa Cuchiero (Vienna University) for the helpful discussions.

\appendix

\section{Technical Result: Swapping Differentiation with Integration}

When deriving the PDE we swapped the operators $\partial_{x}$ with $\bE_{f}$. This essentially requires taking a limit inside an integral, hence we show this is valid in this setting. A similar result was tackled in \cite{LabordereEtAl2015}*{Lemma A2}, although our proof follows similar ideas to the one presented there, our version relaxes some of the conditions on the second derivative.
\begin{lemma}
	\label{Lem:Swap derivative and integral}
	Let Assumptions \ref{Assume:Main SDE Assumptions}, \ref{Assume:Bounded g} and \ref{Assume:Finite Variance} hold. Let $\psi^{T_{1}, \bar{X}_{T_{1}}}$, $\Delta b_{2} \tilde{\psi}^{T_{2},\bar{X}_{T_{2}}} \cW_{3}^{1}$ and
	$\sigma(\theta_{1}, T_{2})^{2}  \tilde{\psi}^{T_{2},\bar{X}_{T_{2}}} \cW_{3}^{2}$ be $\bP$-a.s. uniformly integrable,  let $\tilde{\psi}^{T_{1},\bar{X}_{T_{1}}} \cW_{2}^{1}$ and $\tilde{\psi}^{T_{1},\bar{X}_{T_{1}}} \cW_{2}^{2}$ be $\bP$-a.s. integrable (as defined in Theorem \ref{Thm:Representation Solves the PDE}), and define the function
	\begin{align*}
	\hat{v}(T_{1}, \bar{X}_{T_{1}}) 
	:=
	\bE_{f}[\bE_{W}[ \psi^{T_{1}, \bar{X}_{T_{1}}} | \cF_{T_{1}}] | \cF_{T_{1}}] \, .
	\end{align*}
	Then for $i \in \{1,2\}$,
	\begin{align*}
	\partial_{x}^{i}\hat{v}(T_{1}, \bar{X}_{T_{1}}) 
	=
	\bE_{f}[\partial_{x}^{i}\bE_{W}[ \psi^{T_{1}, \bar{X}_{T_{1}}} | \cF_{T_{1}}] | \cF_{T_{1}}]
	\quad
	\bP \text{-a.s.}
	\end{align*}
\end{lemma}

\begin{proof}
	Technically, the results below are for random variables and hence should be viewed in the a.s. sense, however, for ease of presentation we suppress writing a.s. at the end of each equation.
	Let us start by noting that,
	\begin{align*}
	\psi^{T_{1},\bar{X}_{T_{1}}}
	=
	\1_{\{N_{T}(T_{1}) =0 \}}
	\frac{g(\bar{X}_{T_{N_{T}+1}})}{\overline{F}(\Delta T_{N_{T}+1})}
	+
	\1_{\{N_{T}(T_{1}) \ge 1\}}
	\beta 
	\prod_{k=3}^{N_{T}}
	P_{k} \, ,
	\end{align*}
	where $N_{T}(T_{1})= N_{T} - N_{T_{1}}$. Observing that we can remove the time integral for the case $N_{T}(T_{1})=0$, that is,
	\begin{align*}
	\hat{v}(T_{1}, \bar{X}_{T_{1}}) 
	=
	\bE_{f}[\bE_{W}[ \1_{\{N_{T}(T_{1}) =0 \}}\psi^{T_{1}, \bar{X}_{T_{1}}} 
	+ \1_{\{N_{T}(T_{1}) \ge 1 \}}\psi^{T_{1}, \bar{X}_{T_{1}}}
	| \cF_{T_{1}}] | \cF_{T_{1}}] \, ,
	\end{align*}
	and by integrability we have 
	\begin{align*}
	\bE_{f}[\bE_{W}[ \1_{\{N_{T}(T_{1}) =0 \}}\psi^{T_{1}, \bar{X}_{T_{1}}} 
	| \cF_{T_{1}}] | \cF_{T_{1}}] 
	&
	= 
	\bE_{W}[\bE_{f}[ \1_{\{N_{T}(T_{1}) =0 \}}\psi^{T_{1}, \bar{X}_{T_{1}}} 
	| \cF_{T_{1}}] | \cF_{T_{1}}]
	\\ &	
	= 
	\bE_{W}[g(\bar{X}_{T_{N_{T}+1}}) | \cF_{T_{1}}] \, .
	\end{align*}
	Hence we only need to consider the case $N_{T}(T_{1}) \ge 1$ hence $T_{2}<T$. To make the proof easier we define the function $\varphi$ for $T_{1} < T_{2} < T$ and $\bar{X}_{T_{2}} \in \bR$ as follows,
	\begin{align*}
	\frac{1}{f(\Delta T_{2})} \varphi^{T_{1}, \bar{X}_{T_{1}}}(T_{2}, \bar{X}_{T_{2}})
	=
	\bE[\1_{\{N_{T}(T_{1}) \ge 1 \}}\psi^{T_{1}, \bar{X}_{T_{1}}}
	| \cF_{T_{2}}] \, .
	\end{align*}
	Following the argument as in Theorem \ref{theo:StochasticRepresentation} one can conclude from our uniform integrability assumption that for any $T_{1}< T_{2} < T$, $\varphi^{T_{1}, \bar{X}_{T_{1}}}(T_{2}, \bar{X}_{T_{2}})$ is $\bP$-a.s. continuous in space i.e. w.r.t. $\bar{X}_{T_{2}}$. Further for any fixed $t<T_{1}<T_{2}$, $\varphi$ is bounded in space. To see this one can observe for $T_{2}<T$,
	\begin{align*}
	|\varphi^{T_{1}, \bar{X}_{T_{1}}}(T_{2}, \bar{X}_{T_{2}}) |
	= &
	|f(\Delta T_{2}) \bE[\1_{\{N_{T}(T_{1}) \ge 1 \}}\psi^{T_{1}, \bar{X}_{T_{1}}}
	| \cF_{T_{2}}] | 
	\\
	= & \Big|f(\Delta T_{2})
	\bE 
	\Bigg[
	\Bigg(
	\frac{\Delta g_{T_{N_{T}+1}}}{2\overline{F}(\Delta T_{N_{T}+1})}
	\frac{\Delta b_{N_{T}} \cW_{N_{T}+1}^{1} - \frac{1}{2} \sigma(\theta_{N_{T}-1},T_{N_{T}})^{2} \cW_{N_{T}+1}^{2}}{f(\Delta T_{N_{T}})} 
	\\
	&
	+ \frac{\Delta \hat{g}_{T_{N_{T}+1}}}{2\overline{F}(\Delta T_{N_{T}+1})}
	\frac{-\Delta b_{N_{T}} \cW_{N_{T}+1}^{1} - \frac{1}{2} \sigma(\theta_{N_{T}-1},T_{N_{T}})^{2} \cW_{N_{T}+1}^{2}}{f(\Delta T_{N_{T}})}
	\Bigg)
	\\
	& \times
	\prod_{k=3}^{N_{T}}
	\frac{\Delta b_{k-1} \cW_{k}^{1} - \frac{1}{2} \sigma(\theta_{k-2},T_{k-1})^{2} \cW_{k}^{2}}{f(\Delta T_{k-1})}
	\Big | \cF_{T_{2}}
	\Bigg] \Big|.
	\end{align*} 
	Removing $\cF_{T_{2}}$-measurable terms and noticing that the remaining terms are integrable and $\Delta b_{k-1}< C$ independent of $\bar{X}_{T_{2}}$, we have $\varphi^{T_{1}, \bar{X}_{T_{1}}}(T_{2}, \cdot)$ is bounded in space, as required. Hence we can consider the following bounded Lipschitz approximation to $\varphi$,
	\begin{align*}
	\varphi_{n}^{T_{1}, \bar{X}_{T_{1}}}(T_{2}, x)
	:=
	\inf_{y \in \bR}
	\Big\{
	\varphi^{T_{1}, \bar{X}_{T_{1}}}(T_{2}, y)
	+
	n|x-y|
	\Big\} \, .
	\end{align*}
	One can observe this approximation is both pointwise convergent and increasing in $n$. We therefore work with this approximation and take the limit to complete the proof.
	
	Let us consider differentiating w.r.t. $x$, and in order to make all steps clear let us explicitly write each expectation. Using the tower property to write $\hat{v}$ in terms of $\varphi$ then making the approximation we obtain,
	\begin{align*}
	& \partial_{x} \bE_{f}
	\left[
	\bE_{W}
	\left[
	\1_{\{N_{T}(T_{1}) \ge 1 \}} \frac{1}{f(\Delta T_{2})} \varphi_{n}^{T_{1}, \bar{X}_{T_{1}}}(T_{2}, \bar{X}_{T_{2}})
	\Big| \cF_{T_{1}} 
	\right]
	\Big| \cF_{T_{1}} 
	\right]
	\\
	&
	=
	\lim_{\epsilon \rightarrow 0}
	\bE_{f}
	\left[
	\frac{1}{\epsilon}
	\bE_{W}
	\left[
	\1_{\{N_{T}(T_{1}) \ge 1 \}} \frac{1}{f(\Delta T_{2})} \big(\varphi_{n}^{T_{1}, \bar{X}_{T_{1}}+\epsilon}(T_{2}, \bar{X}_{T_{2}}^{\epsilon})
	-
	\varphi_{n}^{T_{1}, \bar{X}_{T_{1}}}(T_{2}, \bar{X}_{T_{2}})
	\big)
	\Big| \cF_{T_{1}} 
	\right]
	\Big| \cF_{T_{1}} 
	\right]\, ,
	\end{align*} 
	where we are using the notation $\bar{X}_{T_{2}}^{\epsilon}$ to denote the SDE with initial condition perturbed by $\epsilon$. Dominated convergence theorem implies we can take the limit inside the expectation if we show the ``integrand'' to be bounded. Using the Lipschitz assumption on $\varphi_{n}$, one has that,
	\begin{align*}
	|\varphi_{n}^{T_{1}, \bar{X}_{T_{1}}+\epsilon}(T_{2}, \bar{X}_{T_{2}}^{\epsilon})
	-
	\varphi_{n}^{T_{1}, \bar{X}_{T_{1}}}(T_{2}, \bar{X}_{T_{2}}) |
	\le 
	C|\bar{X}_{T_{2}}^{\epsilon} -\bar{X}_{T_{2}}| \, .
	\end{align*}
	As stated in \cite{LabordereEtAl2015}*{Lemma A2}, since $\bar{X}$ has constant coefficients the following bound holds,
	\begin{align}
	\label{Eq:Epsilon perturbation bound}
	\bE \left[
	\Big|
	\frac{\bar{X}_{T_{2}}^{\epsilon} -\bar{X}_{T_{2}}}{\epsilon}
	\Big|^{2}
	~
	\Big |
	\cF_{T_{1}}
	\right]
	\le
	C \, ,
	\end{align}
	further, since $1/ f(\Delta T_{2}) \le C$ by dominated convergence theorem we can take the limit inside $\bE_{f}$ to conclude,
	\begin{align*}
	& \partial_{x} \bE_{f}
	\left[
	\bE_{W}
	\left[
	\1_{\{N_{T}(T_{1}) \ge 1 \}} \frac{\varphi_{n}^{T_{1}, \bar{X}_{T_{1}}}(T_{2}, \bar{X}_{T_{2}})}
	{f(\Delta T_{2})}
	\Big| \cF_{T_{1}} 
	\right]
	\Big| \cF_{T_{1}} 
	\right]
	\\
	&
	\qquad
	=
	\bE_{f}
	\left[
	\partial_{x}
	\bE_{W}
	\left[
	\1_{\{N_{T}(T_{1}) \ge 1 \}} \frac{\varphi_{n}^{T_{1}, \bar{X}_{T_{1}}}(T_{2}, \bar{X}_{T_{2}})}
	{f(\Delta T_{2})}
	\Big| \cF_{T_{1}} 
	\right]
	\Big| \cF_{T_{1}} 
	\right]\, .
	\end{align*}
	Completing the proof for the first derivative requires showing one can take the $\lim_{n \rightarrow \infty}$, however, we suppress this here and concentrate on the second derivative. One can check this holds by following the arguments presented in the case of the second derivative.
	
	Again using the sequence of bounded Lipschitz functions we consider,
	\begin{align*}
	& \partial_{x}^{2} \bE_{f}
	\left[
	\bE_{W}
	\left[
	\1_{\{N_{T}(T_{1}) \ge 1 \}} \frac{\varphi_{n}^{T_{1}, \bar{X}_{T_{1}}}(T_{2}, \bar{X}_{T_{2}})}
	{f(\Delta T_{2})}
	\Big| \cF_{T_{1}} 
	\right]
	\Big| \cF_{T_{1}} 
	\right]
	\\
	&
	\qquad
	=
	\lim_{\epsilon \rightarrow 0}
	\bE_{f}
	\left[
	\frac{1}{\epsilon}
	\bE_{W}
	\left[
	\1_{\{N_{T}(T_{1}) \ge 1 \}} \frac{\varphi_{n}^{T_{1}, \bar{X}_{T_{1}}+\epsilon}(T_{2}, \bar{X}_{T_{2}}^{\epsilon})
		-
		\varphi_{n}^{T_{1}, \bar{X}_{T_{1}}}(T_{2}, \bar{X}_{T_{2}})}
	{f(\Delta T_{2})}
	\cW_{2}^{1}
	\Big| \cF_{T_{1}} 
	\right]
	\Big| \cF_{T_{1}} 
	\right]\, ,
	\end{align*}
	where we have used our first derivative result and the fact that $\varphi_{n}$ is a bounded Lipschitz function to rewrite this derivative with a Malliavin weight. To bound this term one can apply Cauchy-Schwarz, use \eqref{Eq:Epsilon perturbation bound} and,
	\begin{align*}
	\bE_{W}
	\left[
	\Big( \frac{\cW_{2}^{1} }
	{f(\Delta T_{2})}
	\Big)^{2}
	\Big| \cF_{T_{1}} 
	\right]
	\le 
	C \, .
	\end{align*}
	Hence we can again apply dominated convergence theorem to obtain,
	\begin{align*}
	& \partial_{x}^{2} \bE_{f}
	\left[
	\bE_{W}
	\left[
	\1_{\{N_{T}(T_{1}) \ge 1 \}} \frac{\varphi_{n}^{T_{1}, \bar{X}_{T_{1}}}(T_{2}, \bar{X}_{T_{2}})}
	{f(\Delta T_{2})}
	\Big| \cF_{T_{1}} 
	\right]
	\Big| \cF_{T_{1}} 
	\right]
	\\
	&
	\qquad
	=
	\bE_{f}
	\left[
	\partial_{x}^{2}
	\bE_{W}
	\left[
	\1_{\{N_{T}(T_{1}) \ge 1 \}} \frac{\varphi_{n}^{T_{1}, \bar{X}_{T_{1}}}(T_{2}, \bar{X}_{T_{2}})}
	{f(\Delta T_{2})}
	\Big| \cF_{T_{1}} 
	\right]
	\Big| \cF_{T_{1}} 
	\right]\, .
	\end{align*}
To complete the proof we need to also take the $\lim_{n \rightarrow \infty}$, and have the expected values the same. Firstly recall that $\varphi$ is an upperbound for $\varphi_{n}$, hence the result follows from the monotone convergence theorem (see \cite{Williams1991}*{Section 5.3}). Alternatively, one can use the upper bound and uniform integrability results in Theorem \ref{Thm:Integrability Hold} to take the $\lim_{n \rightarrow \infty}$.	 	
\end{proof}



\end{document}